\date{Dec 5, 2000}  
\title {      Values of Brownian intersection exponents II:\\
                       Plane exponents 
}
\author { Gregory F.  Lawler\thanks{Duke University, Research
supported by the National Science Foundation}
\and Oded Schramm\thanks{The Weizmann Institute of Science and
Microsoft Research}
\and Wendelin Werner\thanks{Universit\'e Paris-Sud}
}
\documentclass[naturalnames,12pt]{article}
\usepackage{amsmath}
\usepackage{amsthm}
\usepackage{amsfonts}
\usepackage{graphicx} 
\input labelfig.tex 
\newif\ifhyper\IfFileExists{hyperref.sty}{\hypertrue}{\hyperfalse}
\ifhyper\usepackage{hyperref}\fi

\newif\ifdraft
\drafttrue
\numberwithin{equation}{section}
\numberwithin{figure}{section} 

\newtheorem{theorem}{Theorem}
\numberwithin{theorem}{section}
\newtheorem{corollary}[theorem]{Corollary}
\newtheorem{lemma}[theorem]{Lemma}

\newtheorem{proposition}[theorem]{Proposition}

\newcommand{\R}{\mathbb{R}}
\newcommand{\C}{\mathbb{C}}
\newcommand{\Z}{\mathbb{Z}}
\newcommand{\N}{\mathbb{N}}
\newcommand{\HH}{\mathbb{H}}
\def\Disk{\mathbb{U}}
\def\U{\Disk}
\def\UU{\Disk}
\def\Im{{\rm Im}}
\def\Re{{\rm Re}}
\def \A {{A}}
\def\aa{a}
\def\VV{V}
\def\qq{q}

\def \tx {{\tilde \xi}}
\def \P {{\bf P}}
\def \prob {{\bf P}}
\def \expect {{\bf E}}
\def \p {{\partial}}
\def \E {{\bf E}}

\def\SG/{$SLE_\slepar$}
\def \SLE/{$SLE_6$}
\def\SS{{\mathfrak{K}}}
\def\Btrace{{\mathfrak{B}}}
\def\defeq{:=}
\def \proof {{ \medbreak \noindent {\bf Proof.} }}

\def \lam {\nu}

\def \SLG/ {{\SG}}
\def\zzeta{\delta}
\def\L{{\mathfrak{L}}}
\def \tx {{\tilde \xi}}
\def \F {{\cal F}}

\def\defn#1{{\bf #1}} 
\def\ev#1{{\cal {#1}}} 

\def\slepar{\kappa}
\def\st{:\,}
\def\closure{\overline}
\def\k{\slepar}
\def \bar {\overline}

\def\E{{\bf E}}

\def\tilde{\widetilde}
\def\mobtr/{M\"obius transformation}
\def\hat{\widehat}
\def \l {{\ev G}}

\begin{document}
\maketitle
\begin {abstract}
We derive the exact value 
of intersection exponents between planar Brownian motions
or random walks, confirming predictions from theoretical physics
by Duplantier and  Kwon. 
Let $B$ and $B'$ be independent Brownian motions
(or simple random walks) in the plane, started from distinct points.
We prove that the probability that the paths $B[0,t]$
and $B'[0,t]$ do not intersect decays like
$t^{-5/8}$.  More precisely,
there is a constant $c>0$ such that
if $|B_0 - B_0^\prime| =1$,
for all $t \ge1$,
$$ c^{-1} t^{-5/8} \le  \P
\bigl[ B[0,t] \cap B' [0,t] = \emptyset \bigr] 
\le c  t^{-5/8 }. $$
One consequence is that the set
of cut-points of $B[0,1]$ has  Hausdorff dimension $3/4$
almost surely.
The values of other exponents are also derived.
Using an analyticity result, which is to
be established in a forthcoming paper, 
it follows that the Hausdorff dimension of the outer boundary
of $B[0,1]$ is $4/3$, as conjectured by Mandelbrot.

The proofs are based on a study of \SLE/ (stochastic Loewner evolution
with parameter $6$), a recently discovered process which
conjecturally is the scaling limit of critical percolation
cluster boundaries.
The exponents of \SLE/ are calculated, and they agree
with the physicists' predictions for the exponents
for critical percolation and self-avoiding walks.
{}From the \SLE/ exponents the Brownian intersection
exponents are then derived.
\end {abstract} 

\newpage
\tableofcontents
\newpage

\section {Introduction}

This paper is the follow-up of the paper \cite {LSW1},
in which we derived the exact value of intersection 
exponents between Brownian motions in a half-plane.
In the present paper, we will derive the 
value of intersection exponents
between planar Brownian motions (or simple random walks)
in the whole plane.

This problem is very closely related to the more general 
question of the existence and value of critical exponents 
for a wide class of two-dimensional 
systems from statistical physics,
including
percolation, self-avoiding walks, and other random processes.
Theoretical physics predicts that these systems behave 
in a conformally invariant way in the scaling limit, 
and uses this fact to predict certain critical exponents 
associated to these systems.
We refer to  \cite {LSW1} for a more detailed account on this link
and for more references on this  
subject.

Let us now briefly
describe some of  the results that we shall derive in 
the present paper.
Suppose that $B^1, \ldots, B^n$ are $n\ge 2$ 
independent planar Brownian motions started from 
$n$ different points in the plane.
Then it is easy  to see (using a subadditivity argument)
that there exists  a constant $\zeta_n$ such that 
\begin{equation}  \label{jan27.1}
\P
\bigl[
\forall i \not= j \in \{ 1, \ldots, n \},\ 
B^i [0,t] \cap B^j [0,t] = \emptyset 
\bigr] 
=
t^{-\zeta_n + o (1) }
\end{equation}
when $t \to \infty$.
We shall prove that 
\begin {theorem}
\label {main1}
For all $n \ge 2$,
$$
\zeta_n  = \frac { 4 n^2 -1 } {24} .
$$
\end {theorem}

This result had been conjectured  by Duplantier-Kwon \cite {DK}
(see also, more recently, Duplantier \cite {Dqg}), 
using ideas from theoretical physics (conformal field theory, 
quantum gravity, and analogies with some 
other models for which exponents had
also been conjectured).

It was shown in \cite{BL1} 
that the exponent $\zeta_2$ equals the corresponding
exponent for simple random walks (see also \cite{CM}).
This result was sharpened in \cite{Lwalkcut,
LP}, where estimates were derived up to multiplicative
constants.  It follows from these results and Theorem \ref {main1} 
 that if $S$ and $S'$
denote two independent simple random walks started from 
neighboring vertices in $\Z^2$, then for some constant
$c>0$
$$ 
c^{-1} k^{-5/8} \leq
\P \bigl[S[0,k ] \cap S' [0, k ] = \emptyset \bigr] 
\leq c k^{-5/8},
$$ 
for all $k \ge 1$.
Similarly, \cite{Lcut} and Theorem \ref{main1} imply 
\[ c^{-1} t^{-5/8} \leq
 \P[B^1[0,t] \cap B^2[0,t] = \emptyset] 
      \leq c t^{-5/8},  \]
for all $t \ge 1$, assuming that
the distance between $B^1(0)$ and $B^2(0)$ is $1$, say.

One can define more general exponents, allowing intersection
between some Brownian motions, but forbidding intersection
between different packs of Brownian motions.
For instance,
there exists a constant $\zeta=\zeta ( n,m)$ such that if 
$B^1, \ldots, B^n$ and ${B'}^1, \ldots , {B'}^m$
denote $n+m$ independent planar Brownian motions started
from points such that $B^i (0) \not= {B'}^j (0)$
for all $i\le n$ and $j \le m$, then 
\begin {equation}
\label {zetan,m}
\P \bigl[
\forall i \le n, \forall j \le m ,\ 
B^i [0,t] \cap {B'}^j [0,t] = \emptyset
\bigr]
=
t^{- \zeta + o(1) }
\end {equation}
when $t \to \infty$.
Similarly, one can define exponents $\zeta (n_1, \ldots, n_k)$
corresponding to non-intersection between $k$ packs of 
Brownian motions.

It is easy to see (e.g., \cite {Lmulti})
that there is a natural extention of the definition of $\zeta (n,m)$
to pairs $(n,\lambda)$, 
where $n$ is a positive integer and $\lambda$ is any positive real.
In \cite {LW1}, it is  shown that there
is also a natural definition of 
$\zeta (\lambda_1, \ldots , \lambda_k)$ where the
$\lambda_j$ are positive reals with $\lambda_1,\lambda_2 \ge 1$.
In the present paper, we shall derive the value 
of the exponents for a certain class of
$k$-tuples $(\lambda_1, \ldots , \lambda_k)$
(see Theorem \ref {main3}).
In particular, we shall prove 
\begin {theorem}
\label {main2}
For all real $\lambda \ge 2$,
\begin {equation}
\label {m2}
\zeta (2, \lambda ) =
\frac { \bigl(5 + \sqrt {24 \lambda +1} \bigr)^2 - 4 }{96}\,.
\end {equation}
\end {theorem}

It has been shown by Lawler \cite {Lcut, Lfront, Lmulti, Lbuda} that 
some of these
critical exponents are closely related to the Hausdorff
dimension of exceptional subsets of a planar Brownian path.
Recall that a \defn{cut-point} of a connected set $K$ is the set
of points $x\in K$ such that $K\setminus\{x\}$
is disconnected.
The Hausdorff dimension of the set
of cut-points of the Brownian path $B[0,1]$ is $ 2 - 2 \zeta_2$
almost surely \cite {Lcut}.
Consequently, we get the following corollary from Theorem \ref{main1}.

\begin{corollary}
\label{c-cut}
Let $B$ be Brownian motion in the plane.  Then the Hausdorff
dimension of the set of cut-points of $B[0,1]$ is $3/4$ almost
surely.
\end{corollary}

Recall that the \defn{frontier} of a bounded set $K\subset\R^2$
is its outer-boundary, i.e.,  the boundary of the
unbounded connected component of $\R^2\setminus K$.
The exponents $\zeta ( 2, \lambda)$ are closely related to the 
multifractal spectrum of the Brownian frontier \cite{Lmulti}.
In particular \cite{Lfront},
the Hausdorff dimension $d_F$ of the frontier of $B[0,1]$ is almost surely
$d_F= 2 - \eta_2$,
where $\eta_2 := \lim_{\lambda\searrow 0 } 2\zeta (2, \lambda)$
is called the \defn{disconnection exponent} for two Brownian paths.
(This is not the definition of the disconnection exponent
used in \cite{Lfront}; however, the two definitions are 
equivalent; see \cite{Lstrict, LSWup2}.)
 It had been conjectured by Mandelbrot \cite {M} 
(by analogy with the conjectures for planar
self-avoiding walks) that 
$d_F= 4/3$. 
 Upper and lower bounds for $\eta_2$ from \cite {BL2, Wecp}, combined with 
the fact that $d_F =2 - \eta_2$ \cite {Lfront}, showed that 
$1.015 < d_F < 1.5$. (See also 
\cite {Bal} for another proof of $d_F >1$.)
In the present paper, 
we derive the values of $\zeta (2, \lambda)$ only
for $\lambda \ge 2$, so that we cannot directly apply 
our results to show that $\eta_2 = 2/3$.
However, 
in the subsequent paper 
\cite {LSWan}, we prove that:  
\begin {theorem}\label{allLamb}
The function $\lambda \mapsto \zeta (2, \lambda) $ 
is real analytic on $(0, \infty)$. 
\end {theorem}
 
Combining this with Theorem \ref {main2}
shows immediately that 
(\ref {m2}) holds for all $\lambda > 0$, 
 and therefore $\eta_2 = 2/3$. This
completes the proof of Mandelbrot's conjecture:

\begin {corollary}
\label {4/3}
The Hausdorff dimension of the Brownian frontier is 
almost surely $4/3$.
\end {corollary}

The formula for the multifractal spectrum
of the Brownian frontier also follows.
This formula has been conjectured
in \cite {LW1} as a consequence of 
the conjectures of Duplantier-Kwon \cite {DK} and 
of the functional relations between generalized 
Brownian exponents derived in \cite {LW1}; see also 
recent physics work on this subject by Duplantier \cite {Dqg,Dcm}.

The \defn{pioneer points} of $B$ are the
image under $B$ of the set of times $t$ such that
$B(t)$ is in the frontier of $B[0,t]$.
It has been shown \cite{Lbuda} that
the dimension of the set of pioneer points
of $B$ is $2-\eta_1$, where
$\eta_1:=\lim_{\lambda\searrow 0} 2\zeta(1,\lambda)$.
Below, we show that 
$$
\zeta(1,\lambda) =
 \frac{( 3 +
   \sqrt{24 \lambda + 1} )^2  -4 }{96}
$$
for all sufficiently large $\lambda$ (see (\ref{n=1})).
In \cite{LSWan} it will be proven that $\zeta(1,\lambda)$
is analytic for $\lambda>0$.  Consequently,
by analytic continuation of the above formula
for $\zeta(1,\lambda)$, it follows that $\eta_1=1/4$. 
Hence, using the above mentioned result of \cite{Lbuda},
we obtain

\begin{corollary}\label{pion}
The Hausdorff dimension of the set of pioneer points
of Brownian motion is $7/4$ almost surely.
\end{corollary}

Let us briefly mention that there is a (nonrigorous) 
link  between our results and the conjectures concerning 
two-dimensional self-avoiding walks.
For instance, there is a heuristic argument (see \cite {LW2})
 which uses
the Brownian intersection exponents and explains why the
number of self-avoiding walks
of length $N$ on a planar lattice
increases asymptotically like $N^{11/32} \mu^N$, for some
(lattice-dependent) constant $\mu>1$, as conjectured by
Nienhuis~\cite{Nien} (see also~\cite{MS} for a 
mathematical account).  

\medbreak

Just as in \cite {LSW1}, a central role in the present
paper will be played 
by \SLE/, the stochastic Loewner evolution process with parameter 6,
which is conjectured \cite{S1} to correspond to the scaling
limit of two-dimensional critical percolation cluster boundaries.
In \cite{S1} the processes \SG/ were introduced,
and it was shown that $SLE_2$ is the scaling limit
of loop-erased random walk, assuming the conjecture
that the latter has a conformally invariant scaling limit.

Actually, there are two versions of \SG/.  In the first version,
which we now call \defn{radial} \SG/, one has a set $K_t$ growing
from a boundary point of the unit disk to the interior point $0$,
while in \defn{chordal} \SG/, the set $K_t$
grows from a point in $\R$ to $\infty$ within
the upper half plane.  (The precise definitions
will be recalled in Section \ref{sledef}.)
By applying conformal maps, these processes
can be defined in any simply connected proper 
subdomain of the plane.

After recalling the definition of \SG/,
we study in Section \ref {S.3} some of its properties.
In particular, the \SG/ analogues of the
exponents $\zeta(1,\lambda)$, $\lambda\ge 1$, are computed.  
{}From Cardy's formula for \SLE/
(that we proved in \cite {LSW1}) we then
derive the asymptotic decay of the probability that
\SLE/ crosses a long rectangle without touching the upper
and lower boundaries of this rectangle, and show that chordal
and radial \SLE/ are very closely related.

We then turn our attention to the Brownian intersection 
exponents. In Section \ref {S.5}, the definition
and some properties of the exponents are recalled.
In particular,
it is  explained how to formulate these exponents 
in terms of non-intersection between two-dimensional 
Brownian excursions.
Then, these facts  (properties of \SLE/,
exponents for \SLE/, description of the Brownian exponents
in terms of Brownian excursions, properties
of these Brownian excursions) are combined to derive the 
value of the Brownian intersection exponents.

\section {Preliminaries} \label {Preliminaries}

\subsection {Notation}\label{notation}

Throughout this paper, $\U$ will denote the unit disk  
in the complex plane $\C$.
 $C=C_1= \partial \U$ will denote the unit circle.
For any $r >0$, $C_r = r C_1$
will denote the circle 
of radius $r$ centered at 0.
$\HH$ will denote the upper half-plane $\HH = \{ x+ i y \st
y > 0 \}$. 
When $w \not= w'$ are two points on the unit circle $C$, then 
$a(w,w')$ will denote the counterclockwise arc of $C$ from $w$
to $w'$.

Just as in \cite {LW1,LW2,LSW1}, it will be convenient to 
use $\pi$-extremal distances of quadrilaterals: this just means 
$\pi$ times the usual extremal distance.
(Extremal distance is also known as extremal length.
See, e.g., \cite {A} for the definition and basic properties
of extremal length and conformal maps). 
The $\pi$-extremal distance between sets $A$ 
and $B$ in a set $O$ will be denoted
by $\ell ({A}, {B}, O)$.

Let $f$ and $g$ be functions, and let $l\in\R$ or $l=\infty$.
We say that $ f(x) \sim g(x) $ when $x \to l$,
if $f(x)/ g(x)\to 1$.
We write $f(x) \approx g(x)$, if $\log f(x) / \log g(x)\to 1$,
and we write $ f(x) \asymp g(x)$,
if $f(x) / g(x)$ is bounded above and below by 
positive finite constants when $x$ is sufficiently close to $l$.

\subsection {Radial and chordal SLE}
\label{sledef}

In \cite {LSW1}, we studied chordal \SG/ as a random increasing
family $(K_t, t \ge 0)$ of bounded subsets of the upper
half-plane $\HH$
(or, more generally, their image  under a conformal map).
As $t$ increases, the set $K_t$ grows, and $\bigcup_t K_t$
is unbounded.  One can say that $K_t$ is growing towards
$\infty$, which we think of as a boundary point of $\HH$.

Similarly, when looking at the conformal image of $(K_t, 
t \ge 0)$ under the map $ \varphi (z) =  (z-i)/(z+i)$ that maps
$\HH$ onto the unit disk and $\infty$ to $1$, we get an increasing
family of subsets of the unit disk that is growing
towards $1$.
 
In the present paper, we will mainly use a  variant of this 
process, called \defn{radial \SG/}, where $(K_t, t\ge 0)$ is an increasing
family of subsets of the unit disk that 
grows towards $0$. The main distinction is that $0$ is
an interior point of $\U$, instead of a boundary point.

Suppose that $(\zzeta_t, t \ge 0)$ is a 
continuous function taking values on the unit circle $C_1=\p \U$.
Consider for each $z \in \U$, the solution $g_t =g_t(z)$ of the
ordinary differential equation
\begin{equation}\label{e.loewner}
\p_t g_t=
g_t\frac{\zzeta_t +g_t}{\zzeta_t-g_t}\,,
\qquad t\ge 0\,,
\end{equation}
with $g_0=z$.
This equation (and the 
corresponding equation for $g_t^{-1}$)
was first considered by Loewner (see \cite {Lo}, also
\cite {P1}) and is called 
Loewner's differential equation. For each $z \in \U$, it is well-defined  
up to the time
$\tau_z$ where $\lim_{t \nearrow \tau_z} g_t  = \zzeta_{\tau_z}$,
if there is such a $\tau_z$, and otherwise $\tau_z 
=\infty$.
Let $D_t$ be the set of $z \in \U$ such that $t < \tau_z$
(i.e., the set on which $g_t$ is defined), and $K_t = \U \setminus D_t$. 
It is easy to check that $g_t$ is the unique conformal 
map from $D_t$ onto $\U$ such that $g_t (0)=0$ and 
$g_t'(0)$ is a positive real number. (The notation 
$g'$ refers to differentiation with respect to $z$.)
It is also easy to verify that $g_t'(0) = \exp (t)$
by differentiating both sides of  
(\ref{e.loewner}) with respect to $z$ at $z=0$ and
noting that $g_t(0)=0$.

Let $(B_t, t\ge 0)$, be Brownian motion on
the real line, starting at some point $B_0\in\R$,
and let $\slepar\ge 0$.
Set $\zzeta_t\defeq \exp(i \sqrt {\slepar} B_{t})$
and consider the solution to Loewner's differential equation 
as defined above.
(Note that $\zzeta_t$ is just Brownian motion on $\p\U$
with time scaled by $\slepar$.)
The resulting process was defined in \cite {S1} and
called \SG/ (this acronym stands for stochastic Loewner evolution
with parameter $\slepar$). 
The set $K_t$ is called the {hull} of \SG/, and $(\zzeta_t)_{t \ge 0}$
its driving process.

If $f:D\to\U$ is a conformal map,
then radial \SG/ in $D$ starting from $f$ can be defined as the
composition $g_t\circ f$, where $g_t$ is radial \SG/ in $\U$.
Its hull is $f^{-1}(K_t)$, where $K_t$ is the hull of
$g_t$.
If $\p D$ is sufficiently tame, then $f^{-1}$ extends continuously 
to $\p\U$ and hence $f^{-1}(\zzeta_0)$
is well defined, where $\zzeta$ is the driving parameter of
$g_t$.  
We may then refer to the resulting \SG/ process
in $D$ as \SG/ from $f^{-1}(\zzeta_0)$ to $f^{-1}(0)$ in $D$.

In \cite{S1}, another variant SLE 
was also defined (see also \cite{LSW1}),
which we now call \defn{chordal} \SG/.
Let $\HH$ be the upper half plane,
let $\zzeta_t=\sqrt{\slepar}B_t$,
and consider for each $z \in \HH$, the solution
$\tilde g_t =\tilde g_t(z)$ of the
differential equation
\begin{equation}\label{chordal}
\p_t \tilde g_t=
\frac{-2}{\zzeta_t-\tilde g_t}\,,
\qquad t\ge 0\,,
\end{equation}
with $\tilde g_0(z)=z$.
As before, let $\tilde D_t$ be the set of points $z\in\HH$ for which
(\ref{chordal}) has a solution in some interval $[0,t']$,
$t'>t$, and let $\tilde K_t\defeq \HH\setminus \tilde D_t$. 
Then the process $(\tilde g_t)_{t\ge 0}$
is chordal \SG/ in $\HH$ and $\tilde K_t$ is its hull.
Recall that $\tilde K_t$ is bounded for each $t$, but
$\bigcup_{t\ge 0} \tilde K_t$ is unbounded.
If $f:D\to\U$ is a conformal map, then
$f\circ g_t$ is called chordal \SG/ in $D$.
Its hull is $f^{-1}(K_t)$, where $K_t$ is the
hull of the process $g_t$.
If $\p D$ is sufficiently tame, $f^{-1}(\infty)$
and $f^{-1}(\zzeta_0)$ are well defined.
In this case, we refer to the process $f\circ g_t$ as
\SG/ from $f^{-1}(\zzeta_0)$ to $f^{-1}(\infty)$ in $D$.
This terminology is explained more fully in \cite{LSW1}.

For the main results of this paper, only the case $\slepar=6$
will be used.  As we shall see in Section~\ref{S.equiv},
radial \SLE/ is essentially the same process as chordal \SLE/.
For $\slepar\ne 6$, the radial and chordal \SG/ processes
are  closely related, but the equivalence is weaker.

\section{Annulus crossing exponents for SLE}
\label {S.3}
\subsection{Statement}
\label{s.annstat}

Consider a radial \SG/ process
(with $\kappa \ge 0$) 
in the unit disk $\U$, with driving element $\zzeta_t$ starting
at $\zzeta_0=1$.
As described before, let $g_t$ be the conformal map
$g_t: D_t\to\U$,
with  $g_t(0)=0$ and $g_t'(0)=\exp(t)$.
Let 
$$A_t \defeq \p \U\setminus\closure K_{t}
.$$
It is easy to see that 
 $A_t$ is either an arc on $\p \U$, or $A_t=\emptyset$.

Let $b\ge 0$, and set
\begin{equation}\label{e.lambdadef}
\lam=\lam(\k,b)\defeq
\frac{8b+\k-4+\sqrt{(\k-4)^2+16 b \k}}{16}
\,.
\end{equation}

\begin{theorem}
\label{elannulus}
Let $\slepar> 0$ and $r\in(0,1)$, and let $T(r)$ be the least 
$t>0$ such that $K_t$ intersects the 
circle $C_r$.
Let $ \L(r)$ be the $\pi$-extremal distance
from $C_r$ to $C_1$ in
$\U\setminus K_{T(r)}$
(we set $ \L (r) := \infty$ if $A_t = \emptyset$).
Then for $b\ge 1$,
as $r \to 0$,
$$
\E\Bigl[\exp\bigl( -b \L(r)\bigr)\Bigr]
 \asymp
r^\lam.
$$
\end{theorem}

This theorem gives an analogue of Theorem 3.1 in \cite {LSW1}
for crossing
an annulus.  Its proof and usage will also be similar.

\medbreak 
\noindent
{\bf Remarks.}
The constants implicit in the $\asymp$ notation
may depend on $\kappa$ and $b$.

One can also show that the theorem holds when
$b \in (0,1)$, but we
 only need the case $b\ge 1$ (and $\kappa =6$) in the present paper.
Observe that the case 
$\slepar=0$ is also correct, and easy to verify, for then $K_{T(r)}=[r,1]$.

One should also note 
 that the values of the exponents for $\kappa=2$
fit with the conjecture 
that $SLE_2$ is the scaling limit of two-dimensional loop-erased
random walks \cite {S1}
and 
 the exponents computed by Kenyon \cite {K1,K2,K3}
for loop-erased walks. For instance, 
the definition of loop-erased
random walks suggests that 
 the number of vertices in the loop-erasure of
an $N^2$ step walk is roughly $N^{2-\nu(2,1)}=N^{5/4}$
(loosely speaking, in order for one of the $N^2$ steps
of the simple random walk to remain in the loop-erasure, 
the future of the random walk beyond that step has 
to avoid the past loop-erased walk),
in accordance with Kenyon's results.

Similarly, combining this result (with $\kappa=6$)
 with the non-intersection
exponents between \SLE/ in a rectangle (see the 
discussion at the end of \cite {LSW1}) 
and the restriction property for \SLE/ (which was
proven for chordal \SLE/ in \cite{LSW1} and will
be proven for radial \SLE/ in Section \ref {S.4})
leads to the value of non-intersection exponents between 
independent \SLE/ in an annulus that agree with
predictions for annulus crossings in a  
  critical percolation cluster 
  (see \cite {DS2, Ca3, ADA}).
 
\medbreak
 
The proof of Theorem \ref {elannulus}
will consist of three steps. First, we will obtain
an estimate on $\E\bigl[|g_t' (e^{ix})|^b\bigr]$ for large
 (deterministic) times. We deduce from it a result
concerning the large time
behavior of the arclength of $g_t (A_t)$, 
and then show that this implies Theorem 
\ref {elannulus}.

\subsection{Derivative exponents}

\begin{lemma}\label{l.expoder}
Assume that $\slepar>0$ and $b > 0$.
Let $\ev H(x,t)$ denote the event $\{\exp(ix)\in A_t\}$,
and set
\begin{align}
\nonumber
f(x,t)
&
\defeq
\E\Bigl[\left|g_t'\bigl(\exp(i x)\bigr)\right|^b\,
1_{\ev H(x,t)}\Bigr] \,,\\
\nonumber 
q=q(\k,b)
&
\defeq \frac{\k-4+\sqrt{(\k-4)^2+16 b\k}}{2\k}\,,
\\
\nonumber
h^*(x,t)
&
\defeq
\exp(-t\lam)\bigl(\sin(x/2)\bigr)^q\,,
\end{align}
where $\lam$ is as in~(\ref{e.lambdadef}).
Then there is a constant $c>0$ such that
$$
\forall t\ge 1,  \, \forall x\in(0,2\pi), 
\qquad h^*(x,t)\le f(x,t)\le c\, h^*(x,t)\,.
$$
\end{lemma}

\proof
Let
$\zzeta_t = \exp ( i \sqrt {\slepar} B_t )$
be the driving process of the \SG/, with $B_0=0$. 
For all $x \in (0, 2 \pi)$,
let $Y^x_t$ be the continuous real-valued function
of $t$ which satisfies 
$$ g_t (e^{ix} ) = \zzeta_t \exp ( i Y_t^x ) $$ 
and $Y_0^x=x$.
The function $Y_t^x$ is defined on the set of pairs
$(x,t)$  such that $\ev H(x,t)$ holds.
Since
$g_t$ satisfies Loewner's differential equation
\begin{equation}\label{e.loew}
\p _t g_t(z) = g_t(z) \frac{\zzeta_t+g_t(z)}{\zzeta_t-g_t(z)}\,,
\end{equation}
we find that
\begin{equation}\label{Ysde}
d Y_t^x = \sqrt {\slepar} \  dB_t 
+ \cot(Y_t^x/2)\ dt.
\end{equation}

Let 
$$
\tau^x\defeq 
\inf\Bigl\{t \ge 0 \st Y_t^x\in\{0,2\pi\}\Bigr\}\,
$$
denote the time at which $\exp (ix)$ is absorbed
by $\closure{K_t}$, and define  for all $t < \tau^x$
\begin{equation*}
\Phi_t^x :=
\left| g_t'\bigl(\exp(ix)\bigr)\right|
\,
.\end{equation*}
On $t\ge \tau^x$ set $\Phi_t^x\defeq0$.
Note that on $t<\tau^x$
$$
\Phi_t^x =\p_x Y_t^x .
$$
By differentiating (\ref{e.loew}) with respect to $z$, we find
that for $t<\tau^x$
\begin{equation}\label{e.logder}
\p_t \log \Phi_t^x = -
\frac 1{2\sin^2(Y_t^x/2)}\,,
\end{equation}
and hence (since $\Phi_0^x = 1$),
\begin{equation}\label{Phiis}
(\Phi_t^x)^b =
\exp
\left( - \frac {b}{2} 
\int_0^t \frac {ds} {\sin^{2}(Y_s^x/2)} \right)
\ ,
\end{equation}
for $t<\tau^x$.

We now show that the right hand side of~(\ref{Phiis}) is $0$
when $t=\tau^x$.
For all $x \in (0, 2 \pi)$,
choose $m\in\N$ such that  
$x \in (2^{-m}, 2 \pi - 2^{-m})$, 
and define for all $n \ge m$,
the stopping times
\begin {align*}
\rho_n & \defeq  \inf \bigl\{t>0\st
Y_t^x \notin (2^{-n} , 2 \pi - 2^{-n}) \bigr\}, \\
\rho_n' &\defeq  \inf \bigl\{ t > \rho_n \st
| Y_t^x - Y_{\rho_n}^x | \ge 2^{-n-1} \bigr\}, 
\end {align*}
and the event
$
\ev R_n  \defeq 
\{ \rho_n' - \rho_n > 4^{-n} 
 \}$.
It is easy to see 
(for instance by comparing $Y^x$
with two Bessel processes
and using their scaling property, or alternatively,
by  comparing 
with Brownian motions with constant drifts; see e.g., \cite {RY} for the
definition of Bessel processes)
that 
 there is some constant $c'>0$ such
that for every $ n \ge m$
\begin{equation*}
\P  
[ \ev R_n ]  >c'
\,.
\end{equation*}
The strong Markov property
shows that the events $(\ev R_{m},\ev R_{m+1},\dots)$ are independent, so 
that almost surely, there exist infinitely many values of $n\in\N$
such that $\ev R_n$ holds.
For these values of $n$
$$ \int_{\rho_n}^{\rho_{n+1}}
\frac {dt } {\sin^2 (Y_t^x / 2 ) }
\ge \frac { 4^{-n}} {\sin^2 (2^{-n})}
\ge 
1. $$
Hence,  almost surely,
\begin{equation}\label{hit}
\int_0^{\tau^x} \frac {dt } {\sin^2 (Y_t^x / 2 ) } 
= \infty.
\end{equation}

A similar argument shows that
\begin{equation}\label{sidebd}
\lim_{x\searrow0} f(x,t)=\lim_{x\nearrow2\pi} f(x,t)=0
\end{equation}
holds for all fixed $t>0$.
Suppose, for instance, that 
$x \le 2^{- n_0} \min ( \sqrt {t} , \pi / 2 )$, define 
$ \rho_0''=0$ , $x_0 =x$, and for all $ n\ge 1$,
$$
\rho_n''
= \inf\{ s \ge \rho_{n-1}'' \ : \ 
s = \rho_{n-1}'' + (x_{n-1})^2   
\hbox { or }
| Y^x_s - x_{n-1} | \ge x_{n-1} / 2 \} 
$$
and $x_n = Y_{\rho_n''}^x$.
Clearly, for all $n \le n_0$,
$0 < x_n < \pi /2$ and $\rho_n'' < t$.
By comparing $Y^x$ with Bessel processes or Brownian motions
with constant drift, it is easy to check that for 
some (sufficiently small) $c>0$ (independent of $n$ and $x$), if we 
define
the 
events
${\cal R}_n' := 
\{  \rho_{n}''= \rho_{n-1}'' + (x_{n-1})^2\}$, then 
for
all $  n \le n_0$,
$ 
\P [ {\cal R}_n' | {\cal F}_{n-1}  ] 
\ge
c',
$
where ${\cal F}_{n-1}$
denotes the $\sigma$-field
generated by the events
${\cal R}_1', \ldots , {\cal R}_{n-1}'$.
It therefore easily follows that when $x \searrow 0$,
$\int_0^{\min \{\tau^x,t\}} ds / \sin^2 (Y_s^x / 2)  \to \infty$
in probability, and therefore
also that $f(x,t) \to 0$.

Let $F:[0,2\pi]\to\R$ be a continuous function with
$F(0)=F(2\pi)=0$, which is smooth in $(0,2\pi)$, and set
$$
h(x,t)=h_F(x,t):=
\E\Bigl[(\Phi_t^x)^b\,F(Y_t^x)\Bigr].
$$
By~(\ref{Phiis}) and the general
theory of diffusion Markov processes
(see e.g., \cite {Az}), we know that $h$ is smooth in $(0,2\pi)\times\R_+$.
{}From the Markov property for $Y_t^x$ and~(\ref{Phiis}), it
follows that $h(Y_t^x,t'-t) (\Phi_t^x)^b$ is a local
martingale on $t<\min\{\tau^x,t'\}$.  Consequently,
the drift term of the stochastic
differential $d\bigl(h(Y_t^x,t'-t) (\Phi_t^x)^b\bigr)$
is zero at $t=0$.  By It\^o's formula, this means
\begin{equation}\label{ekol}
\p_t h = \Lambda h\,,
\end{equation}
where
$$
\Lambda h:=
\frac \k 2\, \partial_x^2 h
+\cot(x/2)\, \partial_x h
-\frac b {2 \sin^2 (x/2)}\, h
\,.
$$
It can be verified directly that $h^*$ from the statement of the
lemma solves~(\ref{ekol}).
We therefore {\it choose} 
$$
F(x):= \bigl(\sin(x/2)\bigr)^q, 
$$
and claim that $h^*=h_F$.  Indeed, both satisfy~(\ref{ekol}) on
$[0,2\pi]\times [0,\infty)$,
and $h^*(x,0)=F(x)=h_F(x,0)$ on $[0,2\pi]$. 
Moreover, $F\le 1$ implies that $h_F\le f$ everywhere, 
so $h^*-h_F=0$ on $\{0,2\pi\}\times (0,\infty)$,
by~(\ref{sidebd}).  It is also immediate to verify
that $h_F(x,t)\to 0$ as $(x,t)\to (0,0)$ 
or $(x,t)\to (2\pi,0)$.
Set $M=h_F-h^*$.  Then $M$ is smooth in
$(0,2\pi)\times(0,\infty)$, continuous
on $[0,2\pi]\times[0,\infty)$, and satisfies
$\partial_t M=\Lambda M$.

The proof that $M=0$ can be viewed as  a straightforward application
of the maximum principle.  Fix some $\epsilon>0$,
and suppose that
$M\ge\epsilon$ at some point $(x,t)\in[0,2\pi]\times[0,\infty)$.
Among such points, let $(x_0,t_0)$ be a point with $t_0$
minimal.  It is clear that there must
be such a minimal point and that $x_0\in(0,2\pi)$, $t_0>0$.
At $(x_0,t_0)$ we must have $\partial_t M\ge 0$,
by minimality of $t_0$.  Similarly,
$\partial_x M(x_0,t_0)=0$, $\partial_x^2 M(x_0,t_0)\le 0$
 and $M(x_0,t_0)=\epsilon$.
However, this gives
$0\le\partial_t M(x_0,t_0)=(\Lambda M)(x_0,t_0) \le -
b\epsilon/2\sin^2(x_0/2)$,
by the definition of the operator
$\Lambda$, a contradiction.
Since $\epsilon$ was arbitrary, this gives $M\le0$.
The same argument applied to $-M$ shows that $M\ge0$,
which verifies (the subscript will henceforth be omitted from $h_F$) 
\begin{equation}\label{his}
h(x,t)=h^*(x,t)= \exp(-t\lam)\bigl(\sin(x/2)\bigr)^q
\,.
\end{equation}

As mentioned above, $F\le 1$ implies $h(x,t)\le f(x,t)$. 
Therefore, it remains to prove that for all $t \ge 1$ and $x \in (0, 2
\pi)$,
$ f(x,t) \le c h(x,t)$ for some fixed $c>0$. 
By the  Markov property at time $t-1$, we 
have for $t>1$
$$
h(x,t) = \E\bigl[ ( \Phi_{t-1}^x)^b h(Y_{t-1}^x,1)\bigr],
$$
and similarly with $f$ replacing $h$ on both sides.
Hence, it suffices to prove $c h(x,1)\ge f(x,1)$; that is,
$$
c \E [ ( \Phi_1^x)^b F(Y_1^x)] \ge  \E [ ( \Phi_1^x)^b ].
$$
Let $\sigma_y=\sigma^x_y$ be the first time $s \ge 0$
such that $Y^x_s=y$, and if no such $s$ exists, set $\sigma_y=\infty$.
By (\ref{his}), $h(x,t)$ is a decreasing function of $t$.  
This and the strong Markov property give
\begin{equation}\label{mid}
\E[F(Y^x_1)1_{\{\sigma_y< 1\}} (\Phi_1^x)^b] \ge 
h(y,1) \E[1_{\{\sigma_y< 1\}}(\Phi_{\sigma_y}^x)^b]
.\end{equation}
Let 
$$
\aa:=\min\bigl\{y>0\st y-\cot(y/2)=-2\pi\bigr\}
$$
and consider the event
$$
\ev A\defeq
\bigcap_{s\in[ 0,1 ]} \bigl\{ Y^x_s\in (0,a)\bigr\}.
$$
From~(\ref{Ysde}) it then follows that on the event $\ev A$,
$
\sqrt{\slepar}  B_1 \le \aa - \cot(\aa/2)
= -2\pi$.
Define the Brownian motion $\tilde B$ on $[0,1]$
by $\tilde B_s
 := B_s  - 2s B_1$
and define $\tilde Y^x_s$ by the equation~(\ref{Ysde}),
but with $\tilde B$ replacing $B$.
Note that on $\ev A$ we have
$B_1 <0 $, and hence 
\begin{equation}\label{better}
\forall s\in[0,1]\qquad {\tilde Y}^x_s\ge Y^x_s\,.
\end{equation}
Moreover, given $\ev A$, there is a minimal $s_0\in[0,1]$
with ${\tilde Y}^x_{s_0}=\pi$.
Since~(\ref{better}) holds on $\ev A$,
it follows from~(\ref{Phiis}) that
$$
{\tilde\Phi}_{s_0}^x\ge \Phi^x_{s_0}\ge \Phi^x_1\,,
$$
where $\tilde\Phi$ is the analogue of
$\Phi$ for the process ${\tilde Y}^x$.
Since ${\tilde Y}^x$ has the same law as
$Y^x$, with~(\ref{mid}) this implies
$$
h(x,1)\ge h(\pi,1)
\E[1_{\{\sigma_\pi< 1\}} (\Phi^x_{\sigma_\pi})^b] \ge
h(\pi,1) \E[1_{\ev A} (\Phi^x_1)^b] \,.
$$
The same proof gives this relation with $\ev A$ replaced
by the event $\ev A':=
\bigcap_{s \in [0,1]} \{ Y_s^x \in (2 \pi-a , 2 \pi) \}$.  
However, for the event 
$\ev A'':= \{ \exists s \in [0,1] \st Y_s^x \in [a, 2\pi -a] \}$,
we have
$$
h(x,1) \ge h(\aa,1) \E[1_{\ev A''} (\Phi^x_{1})^b].
$$
Since $\ev A\cup\ev A'\cup\ev A''\supset \ev H(x,1)$ and 
$h(\aa,1)\le
h(\pi,1)$,
we get 
$$
h(x,1)\ge h(\aa,1) f(x,1)/3
.$$
  This completes the proof of the lemma.
\qed

\subsection{Harmonic measure exponents}

By conformal invariance,
 the harmonic measure from $0$ of $A_t$
in $D_t = \U \setminus K_t$ is $L_t/(2\pi)$ where
$L_t$ is the length of the arc $g_t(A_t)$.

\begin{theorem}\label{harmexpo}
Suppose that $\slepar>0$ and $b\ge 1$.
Then, when $t \to \infty$, 
\begin{equation*}
\E[(L_t)^b] 
\asymp
 \exp\bigl(-\lam t\bigr).
\end{equation*}
\end{theorem}

\proof 
We have to relate the behavior of $\left|g_t'(e^{ix})\right|^b$, which
we have studied above, to the behavior of
$$
(L_t)^b=
\left(\int_0^{2\pi}
\left|g_t'(e^{ix})\right|\,dx\right)^b
=
\left(\int_0^{2\pi} \Phi^x_t\,dx\right)^b ,
$$
where we set $\Phi_t^x = 0 $ if $\tau^x \le t$.
By convexity of the function $a\mapsto a^b$, it is clear
that 
\begin{align*}\label{e.haveup}
\frac 1{ (2\pi)^{b}}\E[L_t^b]
&
= \E\left[
\Bigl(\frac 1{2\pi}\int_0^{2\pi} \Phi^x_t\,dx\Bigr)^b 
\right]
\le
\E\left[
\frac 1{2\pi}\int_0^{2\pi}\bigl( \Phi^x_t\bigr)^b\,dx 
\right]
\\&
=
\frac 1{2\pi}\int_0^{2\pi}f(x,t)\,dx 
 \le   c \exp(-\lam t)\,,
\end{align*}
where we have used Lemma~\ref{l.expoder} for the last inequality.
Consequently, we only need to prove
the lower bound for $\E[L_t^b]$.

We will find constants
$c_1, c_2>0$ and an event $\ev U_t^*$ such
that
\[   \E\Bigl[(\Phi_t^\pi)^b \,1_{\ev U_t^*}\Bigr] \geq c_1 e^{-\lam t} ,
\]
and on the event $\ev U_t^*$,
\[   \left| \log \Phi_t^\pi - \log \Phi_t^y \right| \leq
    c_2,\qquad \forall y\in  [\pi ,\pi + c_1] . \]
Then
$$
\E[L_t^b] \ge
\E\left[1_{\ev U_t^*} \Bigl(\int_\pi^{\pi + c_1} \Phi^x_t \, dx
\Bigr)^b\right] 
\ge \Bigl( c_1 e^{-c_2}\Bigr)^b  
\E\bigl[1_{\ev U_t^*}(\Phi_t^\pi)^b\bigr]
     \geq  c_3 e^{-\lam t},
$$
which will prove the theorem.

Assume (with no loss of generality) that $t>3$, and let
$t'$ be the integer in $(t-2,t-1]$.
Define the event
$$   \ev \VV_t = \Bigl\{ \frac{\pi}{2} \leq Y_s^x \leq \frac{3 \pi}{2},
 \;\; \forall s\in[t',t]\Bigr\} . $$
It follows from Lemma~\ref{l.expoder} that
there is some constant $c_4>0$ such
that 
$$
\E\bigl[(\Phi_{t'-1}^x)^b\, 1_{\{Y^x_{t'-1}\in[c_4,2\pi-c_4]\}}\bigr]
\ge c_4 e^{-\lam t}.
$$
This clearly implies
\[   \E\bigl[(\Phi_t^\pi)^b \,1_{\ev \VV_t}\bigr] \geq c_5 e^{-\lam t} ,
\]
for some $c_5>0$.

Let $\ev U_t =\ev  U_t(\alpha)$ be the event
\[  \ev U_t \defeq \bigl\{\alpha e^{-s/8} \leq Y_s^\pi \leq 2 \pi
    - \alpha e^{-s/8} , \;\;\forall s\in[0,t] \bigr\} . \]
We claim that for some $\alpha>0$, and every $t>3$, 
\begin{equation}  \label{jun30.1}
\E[(\Phi_t^\pi)^b \,1_{\ev U_t} \,1_{\ev \VV_t}] \geq \frac{1}{2} c_5
    e^{-\lam t} \,.
\end{equation}
To prove this it suffices to show that for some $\alpha>0$,
\begin{equation}  \label{jun30.2}
\E[
(\Phi_t^x)^b\,1_{\neg\ev U_t} \,1_{\ev \VV_t}] \leq \frac{1}{2} c_5 
    e^{-\lam t}.
\end{equation}
For
$u=0,\ldots,t'-1$,  and
$\alpha \in (0,1/5)$,
let $\ev W_u = \ev W_u(\alpha)$ be the event
\[   \ev W_u \defeq \bigl\{\alpha e^{-u/8} \leq Y_s^\pi \leq
     2 \pi - \alpha e^{-u/8}, \;\;\;\forall s\in[u,u+1]
   \bigr\} . \]
Note that 
$$
\bigcup_{u=0}^{t'-1} \neg\ev W_u \supset \ev \VV_t \cap\neg\ev U_t\,.
$$
Hence,
\begin{equation}\label{e.suma}
 \E\bigl[ (\Phi_t^x)^b\,1_{\neg\ev U_t}\, 1_{\ev \VV_t}\bigr]
 \leq \sum_{u=0}^{t'-1} \E\bigl[(\Phi_t^x)^b \,1_{\neg\ev W_u}\bigr] .
\end{equation}
Note that for $u=0,\ldots,t'-1$,
the strong Markov property shows that,
\[  \E\bigl[(\Phi_t^x)^b  \,1_{\neg\ev W_u} \mid
        \F_{u}\bigr] \leq
    (\Phi_{u}^x)^b f\bigl(\alpha e^{-u/8},t-u-1\bigr)  \]
(here, $\F_u$ denotes the sigma-field generated
by $(Y_s, s \le u)$).
Hence, by Lemma~\ref{l.expoder},
\[  \E\bigl[(\Phi_t^x)^b  \,1_{\neg\ev W_u} \bigr] \leq
c_6 \alpha^{\qq}e^{-u\qq/8}e^{-\lam t} .\]
Now~(\ref{e.suma}) gives
$$
 \E\bigl[ (\Phi_t^x)^b\,1_{\neg\ev U_t}\, 1_{\ev \VV_t}\bigr]
\leq c_7\alpha^{\qq} e^{-\lam t},
$$
and hence by choosing $\alpha$ sufficiently small, we get
(\ref{jun30.2}) and therefore (\ref{jun30.1}).
Fix such an $\alpha\in(0,1/5)$, and
let $\ev U_t^* \defeq \ev U_t \cap \ev \VV_t$.

Observe that~(\ref{Ysde})
implies that if $0<x<y < 2 \pi$,
$$
\p_t (Y_t^y-Y_t^x)
= \cot(Y_t^y/2) - \cot(Y_t^x/2)
\le -(Y_t^y-Y_t^x)/2\,,
$$
(since $\cot'(u)\le -1$ in the range $u\in(0,\pi)$)
as long as $t<\min\{\tau^x,\tau^y\}$, 
so that
\begin{equation}\label{e.contract}
\left|Y_t^y-Y_t^x\right|
\le \left|x-y\right| e^{-t/2}.
\end{equation}
Let $y\in(\pi,\pi+\alpha/2)$.
Then
\[       0 < Y_s^y - Y_s^\pi \leq e^{-s/2}(y-\pi)
\le \alpha e^{-s/2} /2 ,
\;\;\;  \forall s \le \min \{\tau^\pi , \tau^y\} . \]
On the event $\ev U_t^*$,
we must  therefore 
 have $t<\min\{\tau^x,\tau^y\}$
and 
$$Y_s^\pi,Y_s^y \in 
\bigl[\alpha e^{-s/8}/2, 2 \pi - \alpha e^{-s/8}/2\bigr]
$$
for all $s\in[0,t]$.
By~(\ref{e.logder}),
this shows that on the event ${\ev U_t^*}$,
for all $s \le t$,
\begin{align*}
&\left|
\p _s \bigl(\log\Phi_s^y- \log \Phi_s^\pi\bigr)\right|
\\
&\quad\le
\left | Y_s^y-Y_s^\pi\right|
\max\Bigl\{\frac12\left| \p_x\bigl( \sin^{-2}(x/2)\bigr)\right|\st
x\in
\bigl[\alpha e^{-s/8}/2, 2 \pi - \alpha e^{-s/8}/2\bigr]\Bigr\}
\\
&\quad \le
c_7 \left | Y_s^y-Y_s^\pi\right| \,e^{3s/8}.
\end{align*}
Now~(\ref{e.contract}) gives
$$
\left|\p _s\bigl(\log\Phi_s^y- \log \Phi_s^\pi\bigr)\right|\le
c_8 e^{-s/8}.
$$
Therefore, on the event $\ev U_t^*$,
\[
\left|\log\bigl(\Phi_t^y/ \Phi_t^\pi\bigr)\right|
\le 8 c_8 . \]
This completes the proof of the theorem.
\qed

\subsection{Extremal distance exponents}\label{s.el}

{\noindent\bf Proof of Theorem~\ref{elannulus}.}
Let $\rho(t)\defeq \inf\bigl\{|z|\st z\in K_t\bigr\}$.
Recall that 
$$
T(r)=\inf\{t\st \rho(t)=r\},
$$
$A_t=\p\U\setminus \closure{K_t}$ and
that $L_t$ is the length of the arc $g_t(A_t)$.
Recall that the Schwarz Lemma says that if
$G:\U\to\U$ is analytic, then $|G'(0)|\le 1$
and Koebe's $1/4$ Theorem says that if $G:\U\to\C$ is conformal with
$G(0)=0$, then $G(\U)\supset (1/4)\left|G'(0)\right|\U$.  
(See, e.g., \cite {A}.)
Since $g_t'(0)=\exp(t)$, applying the former to 
$z\mapsto g_t\bigl(\rho(t) z\bigr)$ 
and the latter to $g_t^{-1}$ give 
$$
\frac 14 e^{-t} \le \rho(t) \le e^{-t}.
$$
In particular, if we fix the radius $r<1/8$ and define 
the deterministic times
$$ t' = t'(r) \defeq   \log (1/r) ,\qquad
t = t(r) \defeq \log (1/ 8r),$$
then
$$
t < T(r) \le t' \le T(r/4).
$$
and 
$$ 
\rho (t) \ge 2r \ge r \ge \rho (t') \ge r/4 \ge r/8.  
$$
In the following lines, 
 $\ell ( S,S';U)$ will stand  for the $\pi$-extremal 
distance between the sets $S$ and $S'$ in $U$. 
Recall that ${\L}(r) = \ell (C_r,C_1 ; \U
\setminus K_{T(r)})$
and define $l_r = \ell (C_r, C_1 ; \U \setminus K_{t(r)})$.
It follows from the above that 
\begin {equation}
\label {i.ed}
l_r
\le 
{\L} (r) 
\le
l_{r/8}.
\end {equation}
Hence, it will be sufficient
 to study the asymptotic 
behavior of $\E [ \exp ( - b l_r )]$.

Note that $g_{t}:D_t\to \U$ is a conformal map defined
on $D_t \supset 2 r\UU$, 
that $g_{t} (0)= 0$ and that 
$g_{t}' (0) = 1/(8r)$.
Hence, it follows immediately from the
Schwarz Lemma and the Koebe $1/4$ Theorem that 
 for any $r < 1/8$,
$$
2^{-5}\UU \subset g_t ( r\UU  ) \subset (1/2)\UU 
.$$
Hence,
$$
\ell ( C_{2^{-5}} , g_t(A_{t}) ; \U ) 
\ge l_r \ge \ell ( C_{1/2} , g_t(A_{t}) ; \U)
$$
and this implies easily that
$$ 
\exp ( - l_r ) \asymp L_{t(r)}.
$$
Theorem~\ref{elannulus} now follows from Theorem~\ref{harmexpo} 
and (\ref {i.ed}).
\qed

\section {Properties of \SLE/}
\label {S.4}
We now turn our attention towards specific properties of 
\SLE/. 

\subsection {Equivalence of chordal and radial \SLE/}
\label{S.equiv}

The following result shows that chordal \SLE/ and radial \SLE/ 
are nearly the same process.
(When $\slepar\ne 6$, a weaker form of equivalence holds.)
A consequence of the equivalence of radial
and chordal \SLE/ is that radial \SLE/
satisfies a restriction property,
since in \cite{LSW1} a restriction property
for chordal \SLE/ has been established.
The significance of the restriction property to the
Brownian intersection exponents has
been evident since \cite{LW2}.
 
\begin{theorem}\label{diskequiv}
Let $\theta\in\p\U\setminus\{1\}$,
and let $K_t$ be the hull of a radial \SLE/ process in the
unit disk $\U$ with driving process $\zzeta_t$ satisfying
$\zzeta_0=\theta$.
Set 
$$
T\defeq\sup\{t\ge 0\st 1\notin\closure{K_t}\}.
$$
Let $\tilde K_u$ be the hull of a chordal \SLE/ process in $\U$ starting
also at $\theta$ and growing towards $1$, and let 
$$
\tilde T\defeq \sup\{u\ge 0\st 0\notin\tilde K_u\}.
$$
Then, up to a random time change, the process $t\mapsto K_t$
restricted to $[0,T)$ has the same law as the process
$u\mapsto\tilde K_u$ restricted to $[0,\tilde T)$.
\end{theorem}

Note that $T$ (resp.\ $\tilde T$) is the first time where 
$K_t$ (resp.\ $\tilde K_u$) disconnects 0 from 1.

\proof
In order to point out where the assumption $\slepar=6$
is important, we let $(K_t, t \ge 0)$
and $(\tilde K_u, u \ge 0)$
be \SG/ processes, without fixing the value of $\slepar$ for the moment.

Let us first briefly recall (see e.g., \cite {LSW1})
how $\tilde K_u$ is defined.
Let 
 $\psi$ be the \mobtr/ that satisfies
$\psi(\U)=\HH$, $\psi(1)=\infty$, $\psi(-1)=0$, and  $\psi(0)=i$;
that is, 
$$
\psi(z)=i\frac {1+z}{1-z}\,.
$$
Suppose that $u \mapsto \tilde B_u$ is a real-valued 
Brownian motion such that  
$\sqrt {\slepar} \tilde B_0 = \psi (e^{i \theta})$.
For all $z \in \U$, define the function $\tilde g_u = \tilde g_u (z)$
such that $\tilde g_0 (z) = \psi (z)$ and 
$$
\partial_u \tilde g_u = \frac {2} { \tilde g_u - \sqrt {\slepar} \tilde
B_u}.
$$
This function is defined up to the (possibly infinite)
 time $\sigma_z$ where 
$\tilde g_u(z)$ hits $\sqrt {\slepar} \tilde B_u$.
Then, $\tilde K_u$ is defined by $\tilde K_u = \{ z \in \U \ : \ \sigma_z 
\le u\}$, so that $\tilde g_u$ is a conformal map from 
$\U \setminus \tilde K_u$ onto the upper half-plane.
This defines
the process
$( \tilde K_u , u \ge 0)$
(the scaling property of Brownian motion shows that the 
choice of the conformal map $\psi$ only influences the law 
of $(\tilde K_u)_{u \ge 0}$ via a time-change).

We are now going to study the radial \SG/.
Let $g_t:\U\setminus K_t\to\U$ be the conformal map 
normalized by $g_t(0)=0$ and $g'_t(0)>0$.
Recall that
\begin{equation}
\label{e.loewd}
\partial_t g =  g \frac{\zzeta+g}{\zzeta-g}\,, 
\end{equation}
where $\zzeta_t= \exp (i \sqrt {\slepar}  B_t)$,
and $B$ is Brownian motion
on $\R$   with $\exp (i B_0 ) = \theta$.
Let $\psi$ be the \mobtr/ as before,
and define
 \begin {align*}
e_t &\defeq g_t (1),
\\
f_t(z) & \defeq  \psi\bigl(g_t(z)/ e_t \bigr),
\\
\gamma _t&\defeq  \psi\bigl(\zzeta_t/ e_t \bigr).
\end {align*}
These are well defined, as long as $t<T$.
Note that 
$f_t$ is a conformal map from $\U \setminus K_t$
onto the upper half-plane, $f_t (1) = \infty$,
and $\gamma_t \in \R$.
{}From (\ref{e.loewd}) it follows that
$$
\partial_t f =-\frac {(1+\gamma^2)(1+f^2)}{ 2(\gamma-f)}.
$$
Let 
$$\phi_t(z)=a(t)z+b(t)$$
 where
$$
a(0)=1,\qquad \partial_ta=-(1+\gamma^2)a/2
$$ 
and
$$
\qquad b(0)=0,\qquad
\partial_tb=-(1+\gamma^2)a\gamma /2.
$$
Set 
\begin {align*}
h_t
&\defeq  \phi_t\circ f_t\,,\\
\beta _t
&\defeq 
\phi_t\bigl( \gamma (t)\bigr)
.\end {align*}
Then
(and this is the reason for the choice of the 
functions $a$ and $b$)
$$
\partial_t  h =
-(a/2)
\frac {(1+\gamma^2)^2}{ \gamma- f}
=
-\frac {(1+\gamma^2)^2a^2/2}{  \beta - h}\,.
$$
$h_t$ is also a conformal map from $\U \setminus K_t$
onto the upper half-plane with $h_t (1) = \infty$. Note also
that $h_0 (z) = \psi (z)$.
We introduce a new time parameter $u = u(t)$ by
setting
$$
\partial_t u = {(1+\gamma^2)^2a^2/4},\qquad u(0)=0\,.
$$
Then
$$
\frac
{\partial h}{\partial u} = \frac {-2}
{   \beta  - 
h}.$$
Since this is the equation defining the chordal \SLE/ process,
it remains to show that
$ u  \mapsto   \beta _{t(u)} / \sqrt {\slepar}$ 
is Brownian motion (stopped at some random time).
This is a direct but  tedious
application of  It\^o's formula:
$$
d\gamma_t = \frac {(1+\gamma^2) \sqrt {\slepar}}{2}\, dB_t
+ \frac {\gamma (1+\gamma^2)}{2} 
\left( \frac {\slepar}{2} - 1  \right)dt 
$$
and 
$$
d  \beta_t 
=
\frac {  (1+ \gamma^2) a}{2}
\left( \sqrt {\slepar}\, dB_t + ( -3 + \frac {\kappa}{2}) \gamma \,dt 
\right).
$$
This proves the claim, and establishes the theorem.
\qed

\medbreak

Note that when $\slepar \not= 6$, even though 
$u\mapsto  \beta$
is not a local martingale, its law is absolutely continuous with
respect to that of $\sqrt {\kappa}$ times a Brownian motion,
as long as $\gamma$ and $u$ remain bounded.
More precisely:

\begin {proposition}
\label {equivalent}
Let $(K_t, t \ge 0)$, $(\tilde K_u, u \ge 0)$, $T$ and $\tilde T$
be defined just as in Theorem \ref {diskequiv}, except that they
are \SG/ with general $\kappa >0$.
There exist two nondecreasing families 
of stopping times $(T_n , n \ge 1)$ and $(\tilde T_n , n \ge 1)$
such that almost surely,
$T_n \to T$ and $\tilde T_n \to \tilde T$ when $n \to \infty$,
and such that for each $n \ge 1$, the laws of 
$ (K_t , t \in [0, T_n])$ 
and $(\tilde K_u , u \in [0, \tilde T_n])$
are equivalent (in the sense that they have a positive density
with respect to each other) modulo increasing time change.
\end {proposition}

\proof
It suffices to take 
$$ T_n = \min \Bigl\{ n,\, \inf \{ t > 0 \ : \ |\zzeta_t - e_t | < 1/ n \}
\Bigr\}. 
$$
Then, it is easy to see that before $T_n$, $|\gamma|$ remains bounded,
$a$ is bounded away from $0$ (note also that $a \le 1$ always), so that 
$t/u$ is bounded and bounded away from $0$.
Hence, $u( T_n)$ is also bounded (since $T_n \le n$).

It now follows directly from Girsanov's Theorem (see e.g., \cite {RY})
that the law of $\bigl(\beta(u) / \sqrt {\kappa}\bigr)_{ u \le  u(T_n)}$
is equivalent to that of 
 Brownian motion up to some (bounded)
stopping time, and Proposition \ref {equivalent} follows.
\qed

\subsection {The crossing exponent for \SLE/}

In this section, we are going to study the probability 
that chordal \SLE/ started at some point 
 on the left-hand side of a 
rectangle crosses the rectangle from the left to the 
right without touching the upper 
and lower boundaries of the rectangle.
As we shall see, the estimate obtained is a direct consequence of 
Cardy's formula for \SLE/ proved in \cite {LSW1}.

\begin{figure} 
\SetLabels
(.5*.5)$\SS$\\
\R(0.02*.65)$w_1$\\
\R(0.0*.45)$w_1'$\\
\L(1.01*.5)$w_3'$\\
\R(0.02*.31)$w_2$\\
\L(0.99*.31)$w_3$\\
\L(0.99*.65)$w_4$\\
\endSetLabels
\centerline{\AffixLabels{\includegraphics*[height=2in]{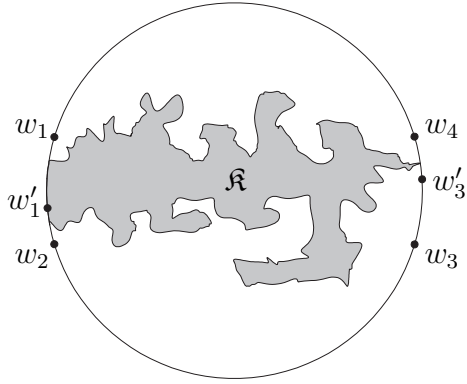}}}
\caption{\label{bang}A succesfull crossing.}
\end{figure}

The notation turns out to be simpler when  considering crossings 
of a quadrilateral in the unit disk, which is equivalent to the 
rectangle case by conformal invariance.
We now describe the setup more precisely.
Recall that when $w,w'\in\p \U$, the couterclockwise arc
from $w$ to $w'$ is denoted $a(w,w')$.
Let $\theta \in (0, \pi/2)$, 
 $\alpha \in (-1, 1)$,
and define the points 
\begin {align*}
w_1 &\defeq \exp ( i ( \pi - \theta ))\,, \\
w_1' &\defeq \exp ( i (\pi + \alpha \theta ))\,, \\
w_2 &\defeq \exp ( i ( \pi + \theta))\,,\\
w_3 &\defeq \exp ( - i \theta )\,, \\
w_4 &\defeq \exp ( i \theta )\,, 
\end {align*}
and note that they appear in counterclockwise order on $\p\U$.
See Figure~\ref{bang}. 
Let $w_3'$ be any point in $a(w_3,w_4)$.
Consider a chordal \SLE/ in $\U$ started from
$w_1'$ and growing towards $w_3'$.
Let $K_t$ be the hull, and let $T$ be the 
first time $t$ such that $\closure{K_t}$
intersects the arc $a(w_3,w_4)$.
Set 
$$
\SS\defeq \bigcup_{t<T}\closure {K_t}
\,.
$$
As shown in \cite {LSW1}, the restriction 
property for \SLE/ shows that up to a monotone time change,
the law of the process $(K_t)_{t<T}$ does
not depend on the choice of $w_3'$.
Since we use the restriction property, 
  the 
result derived in this subsection is specific to $\slepar=6$.

We are interested in the event
$$
\ev E 
\defeq\Bigl \{ {\SS} \cap \bigl(
         a(w_2, w_3) \cup a(w_4, w_1 )\bigr)=\emptyset \Bigr\} \,
.$$

\begin {lemma}
\label {SLEcrossing}
Suppose that $\alpha_0 \in (0, 1)$
is fixed.
When $\theta \searrow 0$, 
$$
\P [ \ev E ]  \asymp  \theta^2 ,
$$
and  this estimate is uniform for $\alpha \in (- \alpha_0, \alpha_0)$.
Moreover,
when $\theta \searrow 0$,
$$
\max_{\alpha \in (-1,1)} \P  [ \ev E ] 
\asymp  \theta^2 .
$$
\end {lemma}

Recall also that, $\theta^2 \sim \exp (- {\ell})$ as
$\theta\searrow0$, where 
${\ell}=\ell\bigl(a(w_1,w_2),a(w_3,w_4);\U\bigr)$
denotes the $\pi$-extremal distance
between $a(w_1, w_2)$ and $a( w_3, w_4)$ in $\U$.

The lemma is hardly surprising.  Suppose for a moment 
that we take $w_1'=-1$ and let $T'$ the first time
such that $\closure{K_t}$ intersects $a(-i,i)$.  Let
$z$ be the point in $a(-i,i)$ which is on the
boundary of $\bigcup_{t<T'} K_t$.  
It is then easy to believe that the probability density
for the location of $z$ should be bounded away
from $0$ and infinity in every compact subset of $a(-i,i)$.
By conformal invariance, this implies the first part of the lemma
for the case $\alpha=0$.  Although it should not be too hard
to prove the lemma with some general arguments such as these,
we find it easier to rely on the more refined results from \cite {LSW1}.
 
\proof
Define the following events:
\begin {align*}
\ev E' &= \bigl\{ {\closure{K_t}} \hbox { hits }
        a(w_2, w_4) \hbox { before } a(w_4, w_1) \bigr\}, \\
\ev E''  &= \bigl\{ {\closure{K_t}} \hbox { hits }
        a(w_2, w_3) \hbox { before } a(w_3, w_1) \bigr\} .
\end {align*}
Note that $ \ev E'' \subset \ev E'$, and that 
$  \ev E = \ev E' \setminus  \ev E''$.
Cardy's formula for \SLE/ derived in \cite {LSW1}
(Theorem 3.2 in the case where $b=0$)
gives the exact value of $\P (\ev E')$ and $\P ( \ev E'')$,
as follows.
Define the cross-ratios 
$$
c' \defeq \frac
{(w_1-w_1')(w_4-w_2)}{(w_4-w_1')(w_1-w_2)}
\,,\qquad
c'' \defeq \frac
{(w_1-w_1')(w_3-w_2)}{(w_3-w_1')(w_1-w_2)}
\,,
$$
and set
$$ 
G ( x) \defeq  
  {}_2F_1( 1/3, 2/3, 4/3 ;\, x )
 \frac {\sqrt{\pi} } {2^{1/3}\Gamma(1/3)\Gamma (7/6) }\, x^{1/3} 
$$
(where ${}_2F_1$ denotes the hypergeometric function). Then
$$ 
\P [\ev E'] = G( c')
\,,\qquad
\P[\ev  E'' ] = G( c'')
\,,\qquad
\P[\ev E]= G(c')-G(c'')
\,.
$$
(To compute $\P [\ev E']$, we view ${K_t}$ as an \SLE/ from $w_1'$ to 
$w_4$, while to compute $\P [\ev E'']$, we view ${K_t}$ as an \SLE/
from $w_1'$ to $w_3$.  As remarked above, the choice of
$w_3'\in a(w_3,w_4)$ does not matter.)
Note that 
$$
c'= \cot\theta \ \tan \bigl((1+\alpha)\theta/2\bigr),\qquad
c'' =\frac {\sin\bigl((1+\alpha)\theta/2\bigr)}
{\sin\theta \cos\bigl((1-\alpha)\theta/2\bigr)}\,.
$$
Both $c'$ and $c''$ converge to 
$ (1+\alpha)/2$ when $\theta \searrow 0$ and
$ c' - c'' \sim (1-\alpha^2)\theta^2  / 4$.
Since $G'(x) =(\sqrt{\pi}/3)\,2^{-1/3}\Gamma(1/3)^{-1}\Gamma(7/6)^{-1}
\bigl((1-x)x\bigr)^{-2/3}$, 
it follows that
$$
\P[\ev E]=
G(c')-G(c'')
\sim 
(\sqrt{\pi}/6)\,\Gamma(1/3)^{-1}\Gamma(7/6)^{-1} 
(1-\alpha^2)^{1/3}\, \theta^2  \,,
$$
as $\theta\searrow 0$, and the lemma follows.
\qed

\section {Brownian intersection exponents}
\label {S.5}

\subsection {Definitions and statement of results}
This section begins with a review of the definitions and some
general facts concerning intersection exponents 
between planar Brownian paths, and proceeds with
a statement of some theorems. 
For more details concerning the
background results on intersection exponents, as well as
some references to the literature,
see \cite {LW1, LW2}.

Suppose that $n+m$
independent Brownian motions 
$B^1, \cdots , B^n$ and $
{B'}^1, \cdots , {B'}^m$ are started from 
points $B^1(0) = \cdots = B^n(0) = i  $ and 
${B'}^1(0) = \cdots = {B'}^m(0) = 1$ in the complex plane,
 and  let ${\Btrace}^j_r, {\Btrace'}^k_r$ denote the traces
(i.e., images) 
of the paths up to the first time they reach the circle $C_r$. 
  Consider
the probability $f_{n,m}(r)$ that the $B$ traces do not
intersect the $B'$ traces, i.e.,
$$ f_{n,m} (r)\defeq 
\prob \left[\Bigl( \bigcup_{j=1}^n \Btrace^j_r \Bigr) \cap
\Bigl(\bigcup_{l=1}^m {\Btrace'}^l_r \Bigr) = \emptyset \right].$$
It is  easy to see that as $r \to \infty$ this probability decays 
like a power law; the $(n,m)$-intersection exponent
$\xi (n,m)$ is defined by
$$ f_{n,m} (r) = r^{- \xi (n,m)  +o(1) },  
\qquad r \to \infty.$$
We call $\xi (n,m)$  the intersection exponent between
one packet of $n$ Brownian motions and one packet of
$m$ Brownian motions.
It is easy to see (e.g., \cite {LW1}) that 
the  exponent $\zeta (n,m)$ described in the introduction
 is equal to 
$\xi (n,m)/2$.

By using the conformal map $z \mapsto 1/z$ and invariance
of planar Brownian motion under conformal mapping, it is clear that 
$$ f_{n,m} (r) = f_{n,m} (1/r)\,,$$
so that the exponents also measure the decay of $f_{n,m}$ 
when $r \searrow 0$.

Similarly, one can define corresponding probabilities
for intersection exponents in a half-plane
$$ \tilde f_{n,m} (r)\defeq 
\prob \left[\Bigl( \bigcup_{j=1}^n \Btrace^j_r \Bigr) \cap
\Bigl(\bigcup_{l=1}^m {\Btrace'}^l_r \Bigr) = \emptyset 
\mbox{ and } 
\Bigl( \bigcup_{j=1}^n \Btrace^j_r \Bigr) \cup
\Bigl(\bigcup_{l=1}^m {\Btrace'}^l_r \Bigr)
 \subset {H}\right],
$$
where ${H}$ is some half-plane through the origin
containing $1$ and $i$.
($\tilde f_{n,m}(r)$ does depend on ${H}$.)
In plain words, we are looking at the probability that all
Brownian motions stay in the half-plane $H$ and that all the $\Btrace$
traces avoid all the $\Btrace'$.
It is also easy to see that there exists a $\tilde \xi (n,m)$ 
(which does not depend on ${H}$) such that 
$$
\tilde f_{n,m} (r) = r^{-\tilde \xi (n,m) +o(1)}
$$
when $r\to\infty$.

These intersection exponents can be generalized in a number
of ways.  For example, 
set 
\[ Z_r :=  \prob\Bigl[{\Btrace'}^1_r \cap 
  \bigcup_{j=1}^n \Btrace^j_r = \emptyset 
   \Bigm|\bigcup_{j=1}^n \Btrace^j_r \Bigr] . \] 
Then $f_{n,m}(r) = \E[Z_r^m]$, and define 
$\xi(n,\lambda)$ for all $\lambda>0$ by the relation 
\[       \E[Z_r^\lambda] \approx r^{-\xi(n,\lambda)} ,\qquad
    r \to \infty . \]
It has been proved by Lawler \cite {Lstrict} that in fact$$ \expect [
Z_r^\lambda ] \asymp r^{-\xi (2, \lambda)}.
$$
The same proof, with minor notational modifications, applies
to the other exponents as well.  
See also \cite {LSWup2} for
a new self-contained proof of this result.

One can also define the exponents $\xi(n_1,n_2,\dots,n_k)$
and $\tx(n_1,n_2,\dots,n_k)$ describing the probability of
non-intersection of $k$ packets of Brownian motions.
In fact, \cite {LW1} proves that there is a natural
and rigorous way to generalize these definitions to the case
where the numbers $n_j$ are positive reals (i.e., not required to be
integers).  For the definition of $\xi(\lambda_1,\dots,\lambda_k)$
one only needs to assume that at least two
of the numbers $\lambda_1,\dots,\lambda_k$ are at least $1$,
and for $\tx$ even this assumption is not necessary.
What makes these extended definitions natural is that they
are uniquely determined by certain identities and relations.
For one, the functions $\xi$ and $\tx$ are invariant under
permutations of their arguments.
Moreover, they satisfy the so-called cascade
relations \cite {LW1}:
for any 
$1 \le  q \le m-1$ and
$(\lambda_1, \ldots , \lambda_m)$ such that 
 $\lambda_1 \ge 1$ and $\max\{\lambda_2,\dots\lambda_m\} \ge 1$,
\begin {equation}
\label {cascade}
\xi ( \lambda_1, \ldots, \lambda_m)
=
\xi \bigl(\lambda_1, \ldots, \lambda_{q}, \tilde \xi (\lambda_{q+1},
\ldots , \lambda_m)\bigr)
. \end {equation}
In \cite {LW1} it was also established that the cascade relations
imply the existence of a continuous
increasing function
$\eta : \bigl[ \tx (1,1) , \infty\bigr) \to \bigl[ \xi (1,1), \infty\bigr)$
such that 
for all $(\lambda_1, \ldots , \lambda_m)\in\R_+^m$ such that at least
two of the $\lambda_j$'s are at least 1,
\begin {equation}
\label {etadef} 
\xi (\lambda_1, \ldots , \lambda_m ) = 
\eta \bigl( \tx ( \lambda_1, \ldots , \lambda_m )\bigr).
\end {equation}

In \cite {LSW1}, we have determined  the exact value of 
exponents $\tx (\lambda_1, \ldots , \lambda_m)$ for a certain   class 
of numbers $(\lambda_1, \ldots, \lambda_m)$, namely, for all
$m \ge 2$, and all
$(\lambda_1, \ldots , \lambda_{m})\in \{ l ( l+1) / 6 \st l \in \N \}^{m-1}
\times \R_+$,
\begin {equation}
\label {txgeneral}
\tx (\lambda_1, \ldots , \lambda_m )
=
\frac {\left( 
\sqrt { 24 \lambda_1 + 1} +
 \cdots + \sqrt { 24 \lambda_m +1 } 
   - (m-1) \right)^2  - 1  }{ 24 } \,.
\end {equation}
In particular
\begin {align}
\label {tx1,1,lambda}
\tx (1 ,  1, \lambda)  
&=
\frac { \left(8 + \sqrt {24\lambda +1 }\right)^2 - 1  } {24}\,,\\
\label {tx1,1,1,lambda}
\tx (1 ,1, 1, \lambda) 
&=
\frac { \left(12 + \sqrt {24\lambda +1 }\right)^2 - 1  } {24}\,.
\end {align}

In the next section, we will
prove the following two results:

\begin {theorem}
\begin {equation}
\label {xi1,1}
\xi (1,1) = 5/4
. \end {equation}
\end {theorem} 

\begin{theorem}   
\label{snowday.theorem}
For every $\lambda \ge 0$,
\begin {equation}
\label {xi1,1,1,lambda}
  \xi(1,1,1,\lambda) 
=
\frac {(11 + \sqrt {24 \lambda +1 })^2 - 4 }{48}
. \end {equation}
\end{theorem}

Let us now see what these two results directly imply:

\begin {theorem}
\label {main3}
For all $x \ge 7$,
\begin {equation}
\label {eta}
\eta (x) = \frac { (\sqrt {24x +1} -1 )^2 - 4 }{48}\,.
\end {equation}
Moreover, for all $m \ge 2$ and for all 
$(\lambda_1, \ldots, \lambda_{m})\in
\{ l(l+1) / 6 \st l \in \N \}^{m-1}\times\R_+$,
\begin {equation}
\label {general}
\xi \left( \lambda_1, \ldots, \lambda_m \right)
=\frac 
{ \left( \sqrt {24 \lambda_1 + 1} + \cdots + \sqrt {24 \lambda_m +1 }-
  m  \right)^2 - 4 
}{ 48}
\end {equation}
provided that at least two of the numbers $\lambda_1, \ldots, \lambda_m$
are at least $1$, and the right-hand side of (\ref {general})
is at least $35/12$.
\end {theorem}

In particular, for all $\lambda \ge 10/3$ and $\lambda' \ge 2$,
\begin {align}
\label {n=1}
\xi (1 , \lambda ) &=
 \frac{( 3 +
   \sqrt{24 \lambda + 1} )^2  -4 }{48}\,,\\
\label {n=2}
\xi (2, \lambda' ) &=
  \frac{(
5+    \sqrt{24 \lambda' + 1} )^2 - 4}{48}\,,
\end {align}
and for all integers $m \ge 2$ (using (\ref{xi1,1}) for the case $m=2$),
\begin {equation}
\label {1otimesN}
\xi (1^{\otimes m}) = \frac {4 m^2 - 1 }{12 }\,,
\end {equation}
where  $1^{\otimes m} \defeq (1, \ldots , 1) \in \N^m$.

\medbreak \noindent {\bf Proof of Theorem \ref{main3}} (assuming
(\ref {xi1,1}) and (\ref {xi1,1,1,lambda})){\bf .}
Combining (\ref {tx1,1,1,lambda}), (\ref {xi1,1,1,lambda})
and (\ref {etadef}) gives (\ref {eta}).
Hence, we get (\ref {general})  from (\ref {txgeneral}), 
(\ref {etadef}) and (\ref {eta}).
\qed

\medbreak
\noindent
{\bf Proof of Theorems \ref{main1} and \ref {main2}}
(assuming (\ref {xi1,1}) and (\ref {xi1,1,1,lambda})){\bf .}
In view of the fact that the time exponents differ by a factor
of 2 from the space exponents, the theorems follow from 
Theorem \ref {main3}.
\qed

\medbreak

\noindent
{\bf Remark.}
In \cite {LSW2s},
we will
show that (\ref {txgeneral}) holds for all 
$(\lambda_1, \ldots , \lambda_m) \in \R_+^m$.
The proofs in the present paper can then be very 
easily   
adapted to  
show that (\ref {general}) holds
for all $(\lambda_1, \ldots , \lambda_m)\in\R_+^m$ such that 
at least two of the $\lambda_j$'s are at least $1$.

\subsection {Excursion measures}
\label{s.excur}

In \cite{LW1,LW2},
a characterization of the intersection exponents was
given in terms of excursions.
The Brownian excursion measures are natural and
interesting objects.  The utility of the excursion
measures in \cite{LW1,LW2} arises from the fact that
in the context of excursion measures, one can study
Brownian motions without specifying the starting
point.  This significantly simplifies some arguments
and estimates.

Roughly speaking, Brownian excursions
in a domain $D$ are Brownian motions started on
the boundary, conditioned to immediately enter $D$,
and stopped upon leaving.  Since Brownian excursions
stay in the domain $D$, conformal transformations of $D$
can be applied to the excursions.  

Let us now describe these Brownian excursions more precisely.
For any bounded simply connected open domain $D$, there exists a
Brownian excursion measure $\mu_D$ in $D$.
This is an infinite measure on the set of paths $(B(t) , t \le \tau)$
in $D$ such that $B(0,\tau) \subset D$ and
$B(0), B(\tau) \in \partial D$ (these endpoints can viewed as prime
ends if necessary).
$x_s\defeq B(0)$ and 
$x_e\defeq B(\tau)$ will denote the starting point and terminal point of the
excursion.
In this discussion, we will 
identify two paths (two excursions) when 
one is obtained by an increasing time-change 
of the other.

One possible definition of $\mu_D$
 is the following. Consider first
$D =\Disk$, the unit disc. For every $s\in(0,1)$ let
$P^s$ be the law of a Brownian motion
started uniformly on the 
circle of radius $s$, and stopped when it exits $\Disk$
(modulo continuous increasing time-change).
Since $z \mapsto \log |z|$ is harmonic, for any $r\in(0,s)$,
$$
P^s[B\hbox{ hits }C_r]=\frac{\log (1/s)}{\log (1/r)}\,.
$$
Set
$$\mu_{\Disk} \defeq 
\lim_{s \nearrow 1} \frac{2 \pi   }{ \log (1/s)}  P^s\, , 
$$
as a weak limit.
Note that the $\mu_{\Disk}$-measure of the set of paths
that hit the circle $C_r$ is $2 \pi / \log (1/r)$.

One can then check 
that for any M\"obius transformation $\phi$ from ${\Disk}
$ onto $\Disk$, 
  $ \phi ( \mu_{\Disk}) = \mu_{\Disk}$.
This makes it possible to extend the definition of $\mu_{D}$
to any 
simply connected domain $D$, by conformal invariance.
These Brownian excursions  
also have a ``restriction''  property 
\cite {LW2}, which is a result of the fact that the Brownian paths
only feel the boundary of $D$ when they hit it (and then stop).

In \cite {LSW1}, we made an extensive use of the Brownian
excursion measure in rectangles ${R}_L = (0, L )
\times (0, \pi)$.
It is easy to see that the measure $\mu_{{R}_L}$,
restricted to those excursions with starting point 
on the left-hand side of the rectangle $[0, i\pi]$,
is obtained as the limit when $s \to 0$
of $\pi s^{-1} P_L^s$, where $P_L^s$ is the law of 
a Brownian motion with uniform starting point on $[s, s + i \pi ]$ which
is stopped when it exits ${R}_L$. 
In particular, this leads to the following result:

\begin {lemma}
\label {lBMcrossing}
Let $\ev E_L$ denote the event that the Brownian excursion $B$ 
in ${R}_L$ 
crosses the rectangle from the left to the right (i.e.,
$x_s \in (0, i \pi )$ and $x_e \in (L, L+i \pi )$). 
Then, when $L \to \infty$,
\begin {equation}
\label {BMcrossing}
\mu_{{R}_L} [\ev E_L ] \asymp e^{-L}.
\end {equation}
\end {lemma}

\proof 
Let $h_z$ denote harmonic
measure from $z$ on $\partial R_L$, where $z\in R_L$.
Since $\Im(\exp z)$ is a harmonic function, it easily
follows that for all $L>1$ and $z\in (1,1+i\pi)$
$$
\sin(\pi/4)\,h_z\bigl([L+i\pi/4,L+3i\pi/4]\bigr)\le  e^{-L}
\le\frac{ h_z\bigl([L,L+i\pi]\bigr)}{\Im(\exp(z))-1}.
$$
It is easy to verify (e.g., by conformal invariance or by a
reflection argument) that there is a constant $c$ such
that $h_z\bigl([L,L+i\pi]\bigr)\le c\,h_z\bigl([L+i\pi/4,L+3i\pi/4]\bigr)$
holds for all $L>2$ and $z\in (1,1+i\pi)$.
Hence, for such $L$ and $z$,
$$
c_1 h_z\bigl([L,L+i\pi]\bigr) \le e^{-L}
\le c_2 h_z\bigl([L,L+i\pi]\bigr) 1_{\{\Im(z)\in[\pi/4,3\pi/4]\}},
$$
with some constants $c_1,c_2>0$.
Since the $\mu_{R_L}$-measure of the excursions which reach
the line $\{\Re(z)=1\}$ does not depend on $L$ and from these
a fixed proportion first hit $\{\Re(z)=1\}$ in $[1+\pi/4,1+3\pi/4]$,
(\ref {BMcrossing}) now follows from the Markov
property, which is valid for the excursion measures.
\qed
 
\bigskip

Similarly, one can define Brownian 
excursions in non-simply 
connected domains. For instance, consider the annulus $\A(r,1)$ 
bounded between the circles
$C_r$ and $C_1$ (where $r\in(0,1)$).
Rather than definining the excursion in $\A(r,1)$ directly,
we base the definition on the excursions in $\U$.
If $\gamma$ is a path, let $\Psi_r(\gamma)$ be
the initial segment of $\gamma$, until the first
hit of $C_r$, or all of $\gamma$,
if $\gamma$ does not hit $C_r$. 
Now set $ \mu_1^r\defeq \Psi_r(\mu_{\Disk})$.
This will be called the Brownian excursion measure on
$\A(r,1)$ for excursions started on $C$.
It is clear that
$\mu_{\Disk}  = \lim_{r \searrow 0} \mu_1^r$.

The measures $\mu_D$ and $\mu_1^r$ are also related
by restriction and conformal invariance.
Suppose that $O$ is a simply connected subset of $\A (r,1)$
such that each of the sets $\bar O \cap C_r$ and $\bar O \cap C$ is
an arc of positive length.  Let 
$L$ denote the $\pi$-extremal distance between these
two arcs in $O$, and let $\phi $ denote the
conformal map from $O$ onto ${R}_L = (0,L) \times (0, \pi)$ 
such that
$\bar O \cap C$ corresponds to $(0,i\pi)$
and
$\bar O \cap C_r$
corresponds to $(L,L+i\pi)$ under $\phi$.

Let $\hat {\ev E}_1$ be the set of paths
starting in $C$ that reach $C_r$ without exiting
$\closure O$.
Consider the image under $\phi$ of the measure $\mu_1^r$
restricted to $\hat{\ev E}_1$.
Then (up to time-change) the image measure 
is exactly $\mu_{{R}_L}$ restricted to the set 
of excursions that cross the rectangle from the left to the right.
(\ref {BMcrossing}) therefore shows that
when $L \to \infty$
\begin {equation}
\label {BEannulus}
\mu_1^r [ \hat {\ev {E}}_1 ] \asymp \exp (-L).
\end {equation}
This may also be easily verified directly.

The mapping  $z \mapsto r/z$ maps 
$\A(r,1)$ conformally onto itself. Let
$\mu_r^1$ be the image of $\mu_1^r$ under this map;
this is a measure on
Brownian excursions started on $C_r$ and stopped
upon leaving $\A(r,1)$.
By symmetry, it follows that
$ \mu_r^1 [ \hat {\ev E}_1 ] =  \mu_1^r [ \hat {\ev E}_1] \asymp \exp (-L)$ 
when $L \to \infty$.

Although we will not use this here, it is worthwhile to note that
the measures $\mu_r^1$ and $\mu_1^r$ agree on the set of paths
crossing the annulus, up to time reversal of the path.

\subsection {Exponents and excursions}
\label{sexpo}
 
We now describe the intersection 
exponents in terms of excursions (referring
to  \cite {LW1,LW2} for the proofs).
Let 
$\bigl(B(t), t \le \tau\bigr)$ be an excursion in ${R}_L$,
and, as above, let
$\ev E_L\defeq \bigl\{\Re(x_s)=0\hbox{ and }\Re(x_e)=L\bigr\}$
be the event that $B$ crosses $R_L$ from left to right.
Let $\Btrace$ denote the image of $B$.
When $\ev E_L$ holds, let $O_B^+$ be the
component of ${R}_L\setminus {\Btrace}$ above ${\Btrace}$,
and let $O_B^-$ be the 
component of ${R}_L \setminus {\Btrace}$ below ${\Btrace}$.
Let ${\L}_B^-$
(respectively ${\L}_B^+$)
denote the $\pi$-extremal distance between $[0, x_s] $ and 
$[L, x_e]$ in $O_B^-$
(respectively $[x_s, i \pi ] $ and $[x_e, L + i \pi ]$
in $O_B^+$).
(We use script fonts for these $\L$ to indicate that they are
random variables.) 
Then, for any $\alpha \ge 0$ and $\alpha' \ge 0$,
the exponent $\tx (\alpha, 1, \alpha') = 
\tx (1, \tx (\alpha, \alpha'))$ is characterized by 
\begin{equation}\label{honey}
\int
_{\ev E_L}
\exp ( - \alpha {\L}^+_B - \alpha' {\L}^-_B ) 
\,d
{\mu_{{R}_L}}(B)
=
\exp \bigl( - \tx (\alpha', 1, \alpha ) L + o(L) \bigr) \,,
\end{equation}
as $L \to \infty$.

Let $\ev E_L^2$ be the set of pairs of
paths $(B,B')\in\ev E_L\times\ev E_L$ such that
the trace $\Btrace'$ of $B'$ is contained in $O_B^-$. 
It follows from (\ref{honey}), 
Lemma \ref {lBMcrossing} and conformal invariance of
the excursion measures that
\begin {equation}
\label {honey2}
\int_{\ev E_L^2} 
\exp ( - \alpha {\L}^+_B )
\,d\mu_{R_L}(B)\,d\mu_{R_L}(B')
= 
\exp \bigl( - \tx (1, 1, \alpha) L + o (L) \bigr)\,,
\end {equation}
as $L \to \infty$.
On $\ev E_L^2$,
let $\L_{B'}^B$ be the $\pi$-extremal distance from $[0,\pi i]$
to $[L,L+\pi i]$ in the domain between $\Btrace$ and $\Btrace'$. 
  Given $B'$, it is clear by
conformal invariance and the restriction property of the
excursion measure that
$1_{\ev E_L^2}\L_{B'}^B$ has the same law as $1_{\ev E_L^2}\L_B^+$.
Consequently, (\ref{honey2}) gives
\begin {equation}
\label {honey3}
\int_{\ev E_L^2} 
\exp ( - \alpha {\L}^B_{B'} )
\,d\mu_{R_L}(B)\,d\mu_{R_L}(B')
= 
\exp \bigl( - \tx (1, 1, \alpha) L + o (L) \bigr)\,. 
\end {equation}

Similarly, one can characterize the exponents $\xi$ in 
terms of the excursion measure $\mu_r^1$.
For any $r<1$, consider two independent 
excursions $B$ and $B'$ of the annulus $\A(r,1)$.
Define the following events:
\begin {align*}
\ev E=
{\ev E}(r) & :=  \{ \Btrace \hbox { crosses the annulus without 
separating } C_r \hbox { from }C \}\,, \\
\tilde {\ev E} =
\tilde {\ev E} (r) & :=  \{ \Btrace \hbox { and } \Btrace' 
\hbox { are disjoint and both cross the annulus} \}\,.
\end {align*}
When ${\ev E}$ holds, let $\L_B$ be the $\pi$-extremal distance
between $C_r$ and $C$ in $\A(r,1) \setminus \Btrace$. 
Similarly, when $\tilde {\ev E}$ is satisfied, let
$O_1$ and $O_2$ be the two components of
$\A(r,1)\setminus (\Btrace\cup \Btrace')$ which have arcs of $C$
on their boundaries, in such a way that
the sequence $\Btrace,O_1,\Btrace',O_2$ is in counterclockwise
order around $\A(r,1)$.
Let $\L_1$ be the $\pi$-extremal
distance from $C_r$ to $C$ in $O_1$, and
let $\L_2$ be the corresponding quantity for $O_2$.
Then for $\lambda,\lambda_1,\lambda_2\in\R_+$ the exponents 
$\xi (1, \lambda)$ and $\xi (1, \lambda_1, 1, \lambda_2)$
can be described as follows \cite {LW2}:
\begin {equation}
\label {BE1}
\int_{\ev E} 
\exp ( - \lambda {\L} )
\,d\mu_r^1(B)
\approx 
r^{\xi (1, \lambda)} ,
\end {equation}
and
\begin {equation}
\label {BE2}
\int_{\tilde {\ev E}}
\exp ( - \lambda_1 {\L}_1 - \lambda_2 {\L}_2 ) \,d\mu_r^1(B)\,d\mu_r^1(B')
\approx r^{\xi (1,\lambda_1, 1, \lambda_2)} 
\,, \end {equation}
as $r\searrow0$.

\subsection {A useful technical lemma}

We now derive a technical
refinement of (\ref {BE1}) and (\ref {BE2})
(in the case $\lambda_1 =1$)
that will be useful to identify the Brownian intersection
exponents with those computed for \SLE/.

Keep the same notation as above, and on ${\ev E}$,
let $\phi$ denote the conformal map that 
maps $O$ onto ${R}_{\L}$ in such a way that the images of
$C_1 \cap \bar O$ and $C_r \cap \bar O$ are mapped onto 
the left and right edges of the rectangle, respectively.
Similarly, on $\tilde {\ev E}$, let $\phi_1$ be
the corresponding conformal map from $O_1$ onto ${R}_{{\L}_1}$.
For all $\alpha>0$ set
\begin {align*}
{\cal H}_\alpha 
&\defeq
{\cal E} \cap \bigl\{  i \in \bar O  \hbox { and }
\phi (i) \in [i \alpha,  i ( \pi - \alpha) ] \bigr\}\,, \\
\tilde {\cal H}_\alpha
&\defeq
\tilde {\cal E} \cap \bigl \{ i \in \bar O_1 
\hbox { and } \phi_1 (i) \in [i \alpha , i (\pi - \alpha)]
\bigr \}
. \end {align*}

\begin {lemma}
\label {refine}
Let $\lambda>0$. 
Then there are sequences $x_n\searrow 0$ and $y_n\searrow 0$
and an $\alpha>0$ such that 
\begin {align}
\label {ref1}
&\int_{{\ev H}_\alpha} \exp ( - \lambda {\L}_B ) \,d\mu_{x_n}^1(B)
\approx 
(x_n)^{ \xi (1, \lambda) } ,\qquad 
&n \to \infty\,,
\\
\label {ref2}
&\int_{\tilde {\ev H}_{\alpha}}
\exp \left( - {\L}_1 - \lambda {\L}_2\right)
\,d\mu_{y_n}^1(B)\,d\mu_{y_n}^1(B')
\approx (y_n)^{\xi (1,1,1,\lambda)} , \qquad
&n \to \infty \,.
\end {align}
\end {lemma}

Actually (see e.g., \cite {LSWup2}), much stronger statements hold:
In the above $\approx$ may be replaced by $\asymp$, 
and these statements hold for every sequence tending to $\infty$.
But the present statement will be sufficient to determine 
the values of the Brownian intersection exponents, and it 
can be easily proved as follows.

\proof
We will only give the detailed proof of (\ref {ref2}). The proof of 
(\ref {ref1}) is easier, follows exactly the same lines,
 and is safely left to the reader.

Because of (\ref {BE2}), it suffices to find the 
lower bound for the left-hand side of (\ref {ref2}).
Let us first introduce some notation.
Let $r\in(0,1/4)$, and consider the measure
$\mu_r^1\times\mu_r^1$ on the space of pairs $(B,B')$.
Let
$$
f(r):=
\int_{\tilde {\ev E}(r)}
\exp \left( - {\L}_1 - \lambda {\L}_2\right)
\,d\mu_{r}^1(B)\,d\mu_{r}^1(B')\,.
$$
Let $B^*$ be the path $B$ stopped when it hits
$C_{1/4}$, if it does, and $B^*=B$, if it
does not hit $C_{1/4}$.
Similarly, define ${B'}^*$ from $B'$.
Let $\Btrace$ and $\Btrace'$ be the traces of $B$ and
$B'$, respectively.
Let $\tilde{\ev E}^*$ be the event that the traces
of $B^*$ and ${B'}^*$ are disjoint.
On $\tilde{\ev E}^*$,
let $O^*_1$ and $O^*_2$ be the domains defined
for $B^*$ and ${B'}^*$, as $O_1$ and $O_2$ were
defined for $B$ and $B'$.
Then $O_j^*\subset O_j$ on $\tilde{\ev E}$, $j=1,2$.
For $j=1,2$, on $\tilde{\ev E}^*$,
let $\L^*_j$ be the $\pi$-extremal distance between
$C_r$ and $C_{1/4}$ in $O_j^*$.
Otherwise, set $\L^*_j=\infty$.

For $a>0$, let $\ev D_a$ be the event
that the distance between $\Btrace\cap \A(1/2,1)$
and $\Btrace'\cap \A(1/2,1)$ is at least $a$.
Suppose that $a\in(0,1/5)$.
Observe that for $(B,B')\notin\ev D_a$,
there is for $j=1$ or $j=2$ a path of length
at most $a$ in
$O_j\cap\A(1/2,1)$, which separates $C_r$ from
$C_1$ in $O_j$.
It then follows that $\L_j\ge\L^*_j+c_1\log (1/a)$,
for some constant $c_1>0$.
Consequently,
\begin{align}
&
\int_{\tilde {\ev E}(r)\setminus \ev D_a}
\exp \left( - {\L}_1 - \lambda {\L}_2\right)
\,d\mu_{r}^1(B)\,d\mu_{r}^1(B')\nonumber
\\&\qquad
\le
a^{c_1\min\{\lambda,1\}}
\int_{\tilde {\ev E}^*}
\exp \left( - {\L}_1^* - \lambda {\L}_2^*\right)
\,d\mu_{r}^1(B)\,d\mu_{r}^1(B')\nonumber
\\&\qquad
=
a^{c_1\min\{\lambda,1\}}
f(4r)\,,
\label{chusham}
\end{align}
since the image of $\mu_r^1$ under the map $B\mapsto 4 B^*$
is $\mu_{4r}^1$.

By (\ref{BE2}), we have
$f(r) \approx r^{\xi (1, 1, \lambda, 1 )}$
(and $f$ is non-decreasing).
Consequently, if $a$ is chosen sufficiently small,
there is a sequence $y_n\searrow 0$ (for instance a subsequence 
of $4^{-n}$) such that
$f(y_n)\ge 2 a^{c_1\min\{\lambda,1\}} f(4y_n)$.
For these $y_n$, (\ref{chusham}) gives
\begin{equation}\label{chubby}
\int_{\tilde {\ev E}(y_n)\cap \ev D_a}
\exp \left( - {\L}_1 - \lambda {\L}_2\right)
\,d\mu_{y_n}^1(B)\,d\mu_{y_n}^1(B')
\ge
a^{c_1\min\{\lambda,1\}}
f(4y_n)\,.
\end{equation}
Fix such an $a$ and such a $y_n$.
Let $\ev I_a$ be the event that $i\in\closure{O_1}$
and the distance from $i$ to $\Btrace\cup\Btrace'$
is at least $a/10$.
Observe that if we apply an independent random uniform rotation
about $0$ to a pair $(B,B')\in\ev D_a\cap\tilde {\ev E}$,
then with probability at least $a/10\pi$ the rotated
pair is in $\ev I_a$.
Since the integrand in (\ref{chubby}) is invariant under rotations,
(\ref{chubby}) and (\ref{BE2}) give
$$
\int_{\tilde {\ev E}(y_n)\cap \ev I_a}
\exp \left( - {\L}_1 - \lambda {\L}_2\right)
\,d\mu_{y_n}^1(B)\,d\mu_{y_n}^1(B')
\approx (y_n)^{\xi (1, 1, \lambda, 1 )}
\,.
$$
It therefore suffices to show that
when $\alpha>0$ is small, we have
$\ev I_a\cap\tilde{\ev E}\subset\tilde{\ev H}_\alpha$.
To prove this, consider a pair $(B,B')\in\ev I_a$.
Let $A$ be the subarc of $\closure{O_1}\cap C_1$
that has $i$ as one endpoint and the other endpoint
is in $\Btrace$, and let $A'$ be the subarc
of $\closure{O_1}\cap C_1$ that has $i$
as one endpoint and the other endpoint
is in $\Btrace'$.  Then the extremal distance from $A$ to
$\Btrace'$ in $O_1$ is bounded from below by a positive constant
depending only on $a$, as is the extremal distance
from $A'$ to $\Btrace$ in $O_1$.
Conformal invariance of extremal distance therefore shows
that the distance from $\phi(i)$ to $\{0,i\pi\}$ is bounded from
below by a constant depending only on $a$,
which proves that $(B,B')\in \tilde{\ev H}_\alpha$,
where $\alpha>0$ depends only on $a$.
\qed

\section {The universality argument} 
\label {S.6}

We are now ready to combine the 
results derived so far to prove our main theorems. 
As in \cite {LSW1}, we 
follow the universality ideas presented 
in \cite {LW2}.
First, Theorem~\ref{xi1,1} ($\xi (1,1) = 5/4$) 
will be proved, followed by
Theorem~\ref{xi1,1,1,lambda} (giving $\xi(1,1,1,\lambda)$).
As we have seen,
Theorems~\ref{main1} and \ref{main2} are immediate consequences.

\subsection {Proof of $\xi (1,1) = 5/4$}

Let $r>0$,
let $K_t$ be a radial \SLE/ process in $\U$ starting at 
$i$, and let $T=T(r)$ be the least $t$ such that
$\closure{K_t}\cap C_r\ne \emptyset$.
Set $\SS:=\closure{K_T}$, let $\P_r$ denote the law
of $\SS$, and let $\expect_r$ denote expectation with respect to
this measure.
As before, we let $\mu_r^1$ denote the Brownian excursion measure 
in $\A(r,1)$ started from $C_r$,
and let $\Btrace$ denote the trace of the excursion $B$.
We are interested in the event ${\ev E}^*$ in which
$\SS\cap \Btrace=\emptyset$ and $\Btrace$ crosses
$\A(r,1)$ (that is, $\Btrace\cap C_1\ne\emptyset$).

The proof will proceed by computing $\bigl(\P_r\times\mu_r^1\bigr)[\ev E^*]$
in two different ways. In the first computation,
we begin by conditioning on $\SS$ and then taking the expectation,
while the second computation begins by conditioning on $B$.

When $\Btrace$ does not separate $C_r$ from $C$ in $\U$,
let $O_B$ denote the connected component of $\A(r,1) \setminus \Btrace$
that touches both circles $C_r$ and $C$
and ${\L}_B$ the corresponding $\pi$-extremal distance.
 
When ${\SS}$ does not disconnect $C_r$ from $C$ in $\U$,
let $O_{\SS}$ denote the connected component 
of $\A (r,1) \setminus {\SS}$ that touches both 
circles and ${\L}_{\SS}$ the corresponding 
$\pi$-extremal distance.
 
Suppose first that ${\SS}$ is given and
that ${\L}_{\SS}<\infty$. 
Note that ${\L}_{\SS}\to\infty$ as $r\searrow 0$.
By the restriction property and conformal invariance of the Brownian 
excursion measure (see Section \ref{s.excur}) and Lemma \ref {lBMcrossing},
$$ 
\mu_r^1 [ B \hbox { crosses the annulus in } O_{{\SS}} ]
= \mu_{R_{{\L}_{\SS}}} [\ev E_{\L_\SS}] 
\asymp  \exp ( - {\L}_{\SS} ),\qquad r\searrow 0\,,
$$
where $\ev E_L$ is as in Lemma \ref{lBMcrossing}.
On the other hand, we know from Theorem \ref {elannulus}
with $b=1,\,\slepar=6$
that 
$$
\expect_r \left[ \exp ( - {\L}_{\SS} ) \right] 
\approx r^{5/4} ,\qquad r\searrow 0\,.
$$
Combining these two facts implies that 
\begin {equation}
\label {estim1}
\bigl(\P_r \times \mu_r^1\bigr) [ {\ev E}^* ] 
\approx r^{5/4} 
,\qquad r\searrow0\,. \end {equation}
 
Suppose now that $B$ is given and that it does
cross the annulus without separating 
the disk $C_r$ from $C$.
The second part of Lemma \ref {SLEcrossing}
 shows that there exists a constant $c>0$ such that 
$$ \P_r [ {\cal E}^* ] \le c \exp ( - {\L}_B ) .$$
Hence, combining this with  (\ref {BE1}) shows that 
\begin{equation}\label{ler}
\bigl(\P_r \times \mu_r^1\bigr)  [ {\ev E}^* ] \le  r^{\xi (1,1) + o(1)}.
\end{equation}

On the other hand, suppose now that $B$ is given
and that $B\in{\ev H}_\alpha$, where ${\ev H}_\alpha$
is as in Lemma~\ref{refine}.
Then by the first part of Lemma \ref {SLEcrossing},
there exists a constant $c'$ such that 
$$
\P_r [ {\ev E}^* ] \ge c' \exp ( - {\L}_B )
.$$
Combining this with Lemma \ref {refine}, we find that
$$
\bigl(\P_{x_n} \times  \mu^1_{x_n}\bigr) 
  [ {\ev E}^* ] \ge
\bigl(\P_{x_n} \times  \mu^1_{x_n}\bigr) 
  [ {\ev H}_\alpha \cap {\ev E}^* ] 
\ge  (x_n)^{\xi (1,1) + o (1)}, \qquad n \to \infty\,.$$
Comparing with (\ref{ler}) gives
$$ \bigl(\P_{x_n} \times  \mu^1_{x_n}\bigr) 
  [  {\ev E}^* ] 
\approx
 (x_n)^{\xi (1,1)}, \qquad n \to \infty\,.$$
Now, by (\ref{estim1}), 
$$ (x_n)^{5/4} 
\approx \bigl(  \P_{x_n} \times  \mu_{x_n}^1\bigr) 
  [ {\cal E}^* ] \approx (x_n)^{\xi (1,1)}$$
when $n \to \infty$, which proves $\xi (1,1) = 5/4$.

\subsection {The determination of $\xi (1,1,1,\lambda)$}

The goal now is to prove (\ref {xi1,1,1,lambda}).
As $\lambda \mapsto \xi (1,1,1,\lambda) = 
\xi (1,1,\tx (1, \lambda))$ 
is continuous on $[0, \infty)$ (see \cite {LW1}), 
we can restrict ourselves to the case where $\lambda>0$.
The proof goes along similar lines as the proof of
$\xi (1,1) = 5/4$. This time, we will 
consider two independent Brownian excursions $B$ 
and $B'$ in the annulus $\A (r,1)$ and one \SLE/ 
${\SS}$ as before.
We are interested in the event that 
$B$, $B'$ and ${\SS}$ all cross the annulus,
that they remain disjoint and that $B$, ${\SS}$ and $B'$ 
are in clockwise order.
Let us call this event $\l$.
Note that in this case, ${\SS}$ crosses the annulus 
in $O_1$, where we use the notations of Section \ref{sexpo}.

We shall compute in two different ways  the quantity
$$
\int_\l
\exp ( - \lambda {\L}_2 ) \,d\mu_r^1(B)\,d\mu_r^1(B')\,d\P_r(\SS)
\,.$$
On the one hand, we know from conformal invariance and the 
restriction property of Brownian excursions and from
(\ref {honey3}) 
 that when ${\SS}$ crosses
the annulus without separating $C_r$ from $C$ 
\begin{equation}\label{eea}
\int_\l
\exp (- \lambda {\L}_2 )\,d\mu_r^1(B)\,d\mu_r^1(B')
\approx
\exp \bigl( - {\L}_{\SS} \tx (1, 1,\lambda ) \bigr)\,,\qquad r\searrow0\,.
\end{equation}
But we know from Theorem \ref {elannulus}
that 
\begin{equation}\label{eeb}
\int \exp ( - b {\L}_{\SS} )\,d\P_r
\approx r^{\nu}
\,,
\end{equation}
where
$$ \nu = \nu ( b) =
\frac { 4b + 1 + \sqrt { 1+ 24b }}{8}\,
.$$
Note that $\nu$ is a continuous increasing function 
of $b$ on $(0, \infty)$.
Consequently, (\ref{eea}) and (\ref{eeb}) give
$$
\int_\l
\exp ( - \lambda {\L}_2 ) \,d\mu_r^1(B)\,d\mu_r^1(B')\,d\P_r(\SS)
\approx r^{ 
\nu (\tx (1,1, \lambda ) ) 
}.$$
Combining this with (\ref {tx1,1,lambda}) shows that 
\begin{equation}
\label{ra}
\int_\l
\exp ( - \lambda {\L}_2 ) \,d\mu_r^1(B)\,d\mu_r^1(B')\,d\P_r(\SS)
\approx r^{
a}  ,
\end{equation}
where
$$
a = a(\lambda):=
 \nu \bigl( \tx(1,1,\lambda)\bigr) =
\frac {\left(11 + \sqrt {24 \lambda +1 }\right)^2 - 4 }{48}
\,.$$

Suppose now that $B$ and $B'$ are given and that 
$\tilde {\cal E}$ is 
satisfied; i.e., that $B$ and $B'$ cross the annulus without 
intersecting each other.
Then Lemma \ref {SLEcrossing}
shows that 
$$\P_r [ {\SS} \subset O_1 ] 
\le c \exp (- {\L}_1 )\,,
$$
for some constant $c$.
Hence, combining this with (\ref {BE2}) shows that 
\begin{equation}\label{abc}
 \int_{\l} \exp (- \lambda {\L}_2 )\,d\P_r(\SS)\,d\mu_r^1(B)
\,d\mu_r^1(B')
\le 
r^{\xi (1,1,1,\lambda) + o (1) }
.\end{equation}

Suppose now that $B$ and $B'$ are given and that 
$(B,B')\in\tilde {\cal H}_{\alpha}$. Then
Lemma \ref {SLEcrossing} shows that there exists a constant $c'>0$ 
such that 
$$ \P_r [ {\SS} \subset O_1  ] 
\ge c' \exp ( - {\L}_1 ). 
$$
Combining this with Lemma \ref {refine}
gives
$$
 \int_{\l} \exp (- \lambda {\L}_2 )\,d\P_{y_n}(\SS)\,d\mu_{y_n}^1(B)
\,d\mu_{y_n}^1(B')
\ge 
(y_n)^{\xi (1,1,1,\lambda) + o (1) }.
$$
Comparing with (\ref{abc}) implies that
$$
 \int_{\l} \exp (- \lambda {\L}_2 )\,d\P_{y_n}(\SS)\,d\mu_{y_n}^1(B)
\,d\mu_{y_n}^1(B')
\approx
(y_n)^{\xi (1,1,1,\lambda) },
$$
when $n \to \infty$.
Consequently, by (\ref{ra}), $\xi (1,1,1, \lambda) = a( \lambda)$.
\qed

\bigskip

\filbreak
\begingroup
\small 
\parindent=0pt 

\ifhyper\def\email#1{\href{mailto:#1}{\texttt{#1}}}\else
\def\email#1{\texttt{#1}}\fi
\vtop{
\hsize=2.3in
Greg Lawler\\
Department of Mathematics\\
Box 90320\\
Duke University\\
Durham NC 27708-0320, USA\\
\email{jose@math.duke.edu}
}
\bigskip
\vtop{
\hsize=2.3in
Oded Schramm\\
Microsoft Corporation,\\
One Microsoft Way,\\
Redmond, WA 98052; USA\\
\email{schramm@microsoft.com}
}
\bigskip
\vtop{
\hsize=2.3in
Wendelin Werner\\
D\'epartement de Math\'ematiques\\
B\^at. 425\\
Universit\'e Paris-Sud\\
91405 ORSAY cedex, France\\
\email{wendelin.werner@math.u-psud.fr}
}
\endgroup

\filbreak

\end {document}